\documentclass[12pt]{article}
\usepackage{amssymb,amsmath,latexsym,amsthm}
\usepackage{url}
\usepackage{enumerate}
\usepackage{ifpdf}
\ifpdf
  \usepackage[pdftex]{hyperref}
\else
  \usepackage[dvips]{hyperref}
\fi

\title{Fluctuations of the empirical quantiles of independent Brownian motions}
\author{Jason Swanson\thanks{Supported in part by NSA Grant H98230-09-1-0079.}\\
University of Central Florida }
\date{August 1, 2010}

\usepackage{color}

\begin{document}

\newtheorem{thm}{Theorem}[section]
\newtheorem{corollary}[thm]{Corollary}
\newtheorem{prop}[thm]{Proposition}
\newtheorem{lemma}[thm]{Lemma}
\theoremstyle{remark}
\newtheorem{remark}[thm]{Remark}
\numberwithin{equation}{section}
\def\al{\alpha}
\def\be{\beta}
\def\ga{\gamma}
\def\de{\delta}
\def\De{\Delta}
\def\ep{\varepsilon}
\def\ze{\zeta}
\def\th{\theta}
\def\ka{\kappa}
\def\la{\lambda}
\def\si{\sigma}
\def\ph{\varphi}
\def\om{\omega}
\def\Om{\Omega}
\def\RR{\mathbb{R}}
\def\NN{\mathbb{N}}
\def\ZZ{\mathbb{Z}}
\def\wt{\widetilde}
\def\wh{\widehat}
\def\ol{\overline}
\def\pa{\partial}
\def\To{\Rightarrow}
\def\pf{\noindent{\bf Proof.} }
\def\qed{\hfill $\Box$}
\providecommand{\flr}[1]{\lfloor#1\rfloor}

\maketitle

\begin{abstract}

We consider iid Brownian motions, $B_j(t)$, where $B_j(0)$ has a rapidly decreasing, smooth density function $f$. The empirical quantiles, or pointwise order statistics, are denoted by $B_{j:n}(t)$, and we consider a sequence $Q_n(t) = B_{j(n):n}(t)$, where $j(n)/n\to\al\in(0,1)$. This sequence converges in probability to $q(t)$, the $\al$-quantile of the law of $B_j(t)$. We first show convergence in law in $C[0,\infty)$ of $F_n=n^{1/2}(Q_n-q)$. We then investigate properties of the limit process $F$, including its local covariance structure, and H\"older-continuity and variations of its sample paths. In particular, we find that $F$ has the same local properties as fBm with Hurst parameter $H=1/4$.

\bigskip


\bigskip


\end{abstract}

\section{Introduction}\label{S:intro}

The classical example of a Brownian motion occurring in nature is the pollen grain on the surface of the water. Bombarded by the much smaller water molecules that surround it, the pollen grain embarks on a random walk with numerous, tiny steps. Applying a law of large numbers type scaling to its path leaves it sitting right where it is, since its trajectory has mean zero. But applying a central limit theorem type scaling, which exposes its fluctuations, reveals the Brownian path that we observe in nature.

If we observe a large collection of pollen grains, and approximate their density with a continuous function, then we might expect this density function to evolve according to the diffusion equation. In this setting, the diffusion equation can be derived by assuming that each individual pollen grain moves as an independent Brownian motion.

However, if the pollen grains are close to one another, as they should be if we are approximating their density with a continuous function, then their motions are certainly not independent. They are interacting with one another through collisions. Moreover, the collisions would evidently provide each individual pollen grain with a drift, pushing it toward regions where the pollen density is lower. In other words, the individual pollen grains would not be performing Brownian motions. The question then arises, what do the trajectories of these colliding pollen grains look like on the mesoscopic scale (that is, under central limit theorem type scaling)?

In this paper, we consider the simpler, one-dimensional question of the scaling limit of colliding Brownian motions on the real line. To motivate our model of the collision process, first consider two physical particles with equal mass and velocities $v_1$ and $v_2$. If these particle have an elastic collision, then they will exchange velocities. The effect of this will be that they exchange trajectories at the moment of collision. If $x_1(t)$ and $x_2(t)$ describe their trajectories without any interaction, then their trajectories in the presence of collisions will be $\min\{x_1(t),x_2(t)\}$ and $\max\{x_1(t),x_2(t)\}$. If we extend this reasoning to $n$ particles, then we are led naturally to consider the empirical quantiles, or order statistics, of the non-interacting trajectories, $x_1(t),\ldots,x_n(t)$. Two concepts, therefore, which are central to our approach, are the notions of an $\al$-quantile of a probability measure, and the order statistics of a family of random variables.

Let $\nu$ be a probability measure on $\RR$ with distribution function $\Phi_\nu (x)=\nu((-\infty,x])$. Given $\al\in(0,1)$, an $\al$-quantile of $\nu$ is any number $q$ such that $\Phi_\nu(q-)\le\al\le\Phi_\nu(q)$. It is easy to see that $\nu$ has at least one $\al$-quantile, given by $q=\inf\{x:\al\le \Phi_\nu(x)\}$. It is also clear that if $\Phi_\nu$ is continuous, then  $\Phi_\nu(q)=\al$ for any $\al$-quantile, $q$.

Given random variables, $X_1,\ldots,X_n$, let $\si$ be a (random) permutation of the set $\{1,\ldots,n\}$ such that $X_{\si(1)} \le\cdots\le X_{\si(n)}$ a.s. For $1\le j\le n$, we define the $j$-th order statistic of $X=(X_1, \ldots, X_n)$ to be $X_{\si(j)}$, and denote it by $X_{j:n}$. Note that  $-X_{j:n} =(-X)_{(n-j+1):n}$.

Now fix $\al\in(0,1)$ and let $B$ be a one-dimensional Brownian motion with a random initial position. For simplicity in the calculations that are to come, we shall assume that $B(0)$ has a density function $f$ that is a Schwartz function. That is, $f\in C^\infty(\RR)$ and
  \begin{equation}\label{Schwartz}
  \sup_{x\in\RR}(1 + |x|^n)|f^{(m)}(x)| < \infty,
  \end{equation}
for all nonnegative integers $n$ and $m$. We will also assume that $f(x)\,dx$ has a unique $\al$-quantile, $q(0)$, such that $f(q(0))>0$.

Let $\{B_j\}$ be an iid sequence of copies of $B$. For fixed $n$, the trajectories $B_1,\ldots,B_n$ denote the paths of a system of $n$ particles with no interactions. When these particles are allowed to interact through collisions, their trajectories are given by the pointwise order statistics, or empirical quantile processes, $B_{1:n}, \ldots, B_{n:n}$. More specifically, $B_{j:n}$ is the process such that, for each $t\ge0$, $B_{j:n}(t)$ is the $j$-th order statistic of $(B_1(t), \ldots,B_n(t))$. Note that $B_{j:n}$ is a continuous process. We shall fix a sequence of integers $\{j(n)\}_{n=1}^\infty$ such that $1\le j(n)\le n$ and $j(n)/n=\al+o(n^{-1/2})$. We then consider the sequence of empirical quantile processes $\{Q_n\}$ given by $Q_n=B_{j(n):n}$. Our first result concerns the convergence in $C[0,\infty)$ of the sequence $\{Q_n\}$. To this end, we begin by defining the (deterministic) limit process, and show that it is continuous.

Let $u(x,t)$ denote the density of $B(t)$. It is well-known that our assumptions on $f$ are sufficient to ensure that $u\in C^\infty (\RR\times(0,\infty))$, $\pa_x^j u\in C(\RR\times[0,\infty))$, and $\pa_x^j u(x,0) =f^{(j)}(x)$ for all $j\ge0$. Moreover, $\pa_x^j u(x,t)=E^x[f^{(j)}(B(t))]$ and $u$ satisfies the diffusion equation, $\pa_t u = (1/2)\pa_x^2 u$. For each $t>0$, the function $u(\cdot,t)$ is strictly positive, which implies that $u(x,t)\,dx$ has a unique $\al$-quantile, $q(t)$. 
The following lemma is easily derived by differentiating the defining equation for $q(t)$; its proof is given in the appendix.

\begin{lemma}\label{L:quant_ODE}
The function $q$ is in $C[0,\infty)\cap C^\infty (0, \infty)$ and satisfies
  \begin{equation}\label{quant_ODE}
  q'(t) = -\frac{\pa_x u(q(t),t)}{2u(q(t),t)}
  \end{equation}
for all $t>0$.
\end{lemma}

\begin{remark}
Let $C^k[0,\infty)$ denote the space of functions $g:[0,\infty)\to \RR$ such that $g^{(j)}$ has a continuous extension to $[0,\infty)$, for all $j\le k$. Also, let $C^\infty[0,\infty)=\cap_{k\ge0}C^k[0, \infty)$. It is easy to see from Lemma \ref{L:quant_ODE} that since $u(q(t),t)>0$ for all $t\ge0$, and $\pa_x^j u\in C (\RR\times[0, \infty))$, it follows that $q\in C^\infty[0,\infty)$. That is, $q^{(k)}$ has a continuous extension to $[0,\infty)$ for all $k$.
\end{remark}

\begin{remark}
The diffusive flux at time $t$, in the positive spatial direction, across an arbitrary moving boundary $s(t)$ is given by $-(1/2)\pa_x u(s(t),t)-u(s(t),t)s'(t)$. The first term corresponds to Fick's first law of diffusion, and the second term comes from the motion of the boundary. By \eqref{quant_ODE}, it follows that $q(t)$ is the unique trajectory starting at $q(0)$, across which there is no diffusive flux. The vanishing of the flux can also be seen by noting that
  \[
  \int_{-\infty}^{q(t)}u(x,t)\,dx = \al,
  \]
which follows from the definition of $q(t)$.
\end{remark}

It will follow from later results that $Q_n\to q$ in probability in $C[0,\infty)$. Our primary interest, however, is with the fluctuations, $F_n=n^{1/2}(Q_n-q)$. Our objective in this paper is twofold: (i) to establish the convergence in law in the space $C[0,\infty)$ of the sequence $F_n$ to a limit process $F$, and (ii) to investigate properties of the limit process $F$, including its local covariance structure, and the H\"older continuity and variations of its sample paths.

Regarding (i), we establish the following result, whose proof can be found in Section \ref{S:param_est}. (The notation $X_n\To X$ means $X_n\to X$ in law.)

\begin{thm}\label{T:main}
There exists a continuous, centered Gaussian process $F$ with covariance function
  \begin{equation}\label{rho_def}
  \rho(s,t) = \frac{P(B(s) \le q(s), B(t) \le q(t)) - \al^2}
    {u(q(s),s)u(q(t),t)},
  \end{equation}
such that $F_n\To F$ in $C[0,\infty)$.
\end{thm}

It is worth pointing out here that the limiting process $F$ is not deterministic at time $t=0$. In fact, $E|F(0)|^2 = (\al - \al^2) f(q(0))^{-2}$.

Regarding (ii), we derive, in Section \ref{S:limit_props}, several key properties of the limit process $F$. We will show that, on compact time intervals, $E|F(t)-F(s)|^2$ is bounded above and below by constant multiples of $|t-s|^{1/2}$. In particular, the paths of $F$ are almost surely locally H\"older continuous with exponent $\ga$ for all $\ga<1/4$. (See Corollaries \ref{C:var_order_0} and Corollary \ref{C:Holder}.) We will show that $F$ is locally anti-persistent. More specifically, nearby, small increments of size $\De t$ have a negative 
covariance
 whose order of magnitude is $-|t-s|^{-3/2}\De t^2$. (See Corollaries \ref{C:IV.1} and \ref{C:IV.2}.) We will also show that $F$ has a nontrivial, deterministic quartic variation, given by $(6/\pi)\int_0^t |u(q(s),s)|^{-2}\,ds$. (See Theorem \ref{T:quart_var}.) All of these are local properties that $F$ shares with $B^{1/4}$, the fractional Brownian motion (fBm) with Hurst parameter $H=1/4$. (Recall that fBm, $B^H$, is a centered Gaussian process with $B^H(0)=0$ and $E|B^H(t) -B^H(s)|^2=|t-s|^{2H}$, where $H\in(0,1)$.) These local properties are also shared with the solution to the one-dimensional stochastic heat equation, which was recently studied in \cite{Sw07.1} and \cite{BS}. On the other hand, the global properties of $F$ can be quite different from $B^{1/4}$, which we will illustrate at the end of Section \ref{S:limit_props}, via the special case where $f(x)\,dx$ is a standard Gaussian distribution.

A model similar to this was studied in \cite{Ha}. That model consists of a countably infinite collection of real-valued stochastic processes $x_i(t)$, $i\in\ZZ$, $t\ge0$. The points $\{x_i(0)\}_{i\in\ZZ}$ form a unit Poisson process on the real line, conditioned to have a point at 0, and labeled so that
  \[
  \cdots < x_{-2}(0) < x_{-1}(0) < x_0(0)
    = 0 < x_1(0) < x_2(0) < \cdots
  \]
The processes $\{x_i(\cdot)-x_i(0)\}$ are independent, standard Brownian motions. This family of processes represents the motion of the particles without collisions. By counting upcrossings and downcrossings, the motion of a ``tagged particle" in the collision system can be defined. Informally, the tagged particle $y(t)$ is defined as follows. It begins with $y(0)=x_0(0)$, and then continues as $y(t)=x_0(t)$ until the first time the path of $x_0$ intersects one of its neighbors, $x_1$ or $x_{-1}$. At this point, $y$ adopts the trajectory of the neighboring particle, and follows this trajectory until the first time it meets one if its two neighbors, and so on.

Of course, when two Brownian particles meet, they intersect infinitely often immediately, which makes it difficult to carry out the above informal construction. This is why upcrossings and downcrossings are used instead. (See also \cite{DGL85}.) Note that, with only finitely many particles, the empirical quantiles serve as tagged particles, without any need for counting crossings.

By the scaling property of Brownian motion, the process $n^{-1}y(n^2t)$ has the same law as a tagged particle in a system initially distributed according to a Poisson process with a density of $n$ particles per unit length. This is analogous to our centered empirical quantile process $Q_n(t)-q(t)$. Multiplying by $n^{1/2}$ shows that our process $F_n$ is analogous to $y_A(t)=y(At)/A^{1/4}$, where $A=n^2$. In \cite{Ha}, it is shown that $y_A(1)=y(A)/A^{1/4}\To(2/\pi)^{1/4}N$ as $A\to \infty$, where $N$ is a standard normal random variable. This implies that for fixed $t>0$,
  \[
  y_A(t) = t^{1/4}y_{At}(1) \To t^{1/4}(2/\pi)^{1/4}N
    \overset{\mathcal{L}}{=} (2/\pi)^{1/4}B^{1/4}(t),
  \]
where $\overset{\mathcal{L}}{=}$ denotes equality in law. In \cite{DGL85}, a substantially stronger result was proven, one case of which shows that $y_A \To(2/\pi)^{1/4}B^{1/4}$ in $C[0,\infty)$. The proof of tightness in \cite{DGL85} relies heavily on the special properties of the initial Poisson distribution, under which the resulting particle system has a stationary distribution.

Similar results hold for the simple symmetric exclusion process (SSEP) on the integer lattice. It has been known since the work of Arratia \cite{A} and Rost and Vares \cite{RVar} in the 1980s that the fluctuations of a tagged particle in a one-dimensional SSEP started in equilibrium converge, in the sense of finite-dimensional distributions
(fdd),
to fBm $B^{1/4}$. The proof of tightness, however, was established only recently by Peligrad and Sethuraman \cite{PS} in 2008. As with the proof in \cite{DGL85}, the proof of tightness for the tagged particle in \cite{PS} also relies on the initial distribution of the particles and does not extend to the non-equilibrium case.

The motion of a tagged particle in a non-equilibrium SSEP was studied in Jara and Landim \cite{JL}. Theorem 2.6 in \cite{JL} establishes a result on the fluctuations of a tagged particle in SSEP which is analogous to our Theorem \ref{T:main}, except that in \cite{JL} convergence is established only in the sense of 
fdds.

It should be noted that the proof methods in \cite{DGL85} do generalize to the non-equilibrium case when one considers the current process, rather than the tagged particle. In 2005, Sepp\"al\"ainen \cite{Se} studied the current process in a system of independent random walks in a non-equilibrium case. Tightness for this current process (for a certain restricted class of initial profiles) was proved by Kumar \cite{Kum} in 2008 by extending the proof of Proposition 5.7 in \cite{DGL85}.

In \cite{Sw07.2}, we considered what is very nearly a special case of the present model. Suppose $B(0)=0$ a.s. and $j(n)=\flr{{(n+1)}/2}$, so that $Q_n$ is the median process, $\al=1/2$, and $q=0$. In that case, it was shown in \cite{Sw07.2} that $F_n\To X$ in $C[0,\infty)$, where $X$ is a centered Gaussian process with covariance function
  \begin{equation}\label{med_covar}
  E[X(s)X(t)]
    = \sqrt{st}\sin^{-1}\left(\frac{s\wedge t}{\sqrt{st}}\right).
  \end{equation}
The main difficulty in establishing this result was the proof of tightness. Unlike \cite{DGL85}, we did not have the benefit of the Poisson initial distribution. Instead, we proved tightness by making direct estimates on $P(|Q_n(t)-Q_n(s)|\ge n^{-1/2}\ep)$. The method of estimation was essentially a four-step procedure. First, we established a connection between this probability and a certain random walk. Second, we derived estimates for probabilities associated with this random walk -- estimates which depend on the parameters of the walk. Third, we estimated those parameters in terms of the dynamics of the original model. And fourth, we put all of this together to get separate estimates on $P(|Q_n(t)-Q_n(s)|\ge n^{-1/2}\ep)$, depending on whether $n^{-1/2} \ep\ll|t-s|^{1/2}$, $n^{-1/2} \ep\approx|t-s|^{1/2}$, or $n^{-1/2} \ep\gg |t-s|^{1/2}$.

Regarding the convergence portion of the present paper, our proof of tightness is motivated by this rough outline. It is our hope that the new ideas and techniques developed here in relation to tightness, as well as the methods employed to study the properties of the limit process, will be applicable to more general colliding particle models.

\section{Finite-dimensional distributions}\label{S:fdd}

We summarize here the definitions and assumptions from Section \ref{S:intro}. Let $B$ be a Brownian motion such that $B(0)$ has a density function $f\in C^\infty$ satisfying \eqref{Schwartz} for all nonnegative integers $n$ and $m$. Let $\{j(n)\}_{n=1}^\infty$ be a sequence of integers satisfying $1\le j(n)\le n$ and $j(n)/n=\al+o(n^{-1/2})$, where $\al\in(0,1)$. Assume the probability measure $f(x)\,dx$ has a unique $\al$-quantile, $q(0)$, satisfying $f(q(0))>0$.

Let $\{B_j\}$ be a sequence of iid copies of $B$. We define $F_n=n^{1/2} (Q_n-q)$, where $Q_n=B_{j(n):n}$ and $q(t)$ is the unique $\al$-quantile of the law of $B(t)$. The fact that $F_n\in C[0,\infty)$ follows from Lemma \ref{L:quant_ODE}. 

In this section, we begin the proof of Theorem \ref{T:main} by noting that the convergence of the fdds follows as an immediate corollary of a multi-dimensional quantile central limit theorem, which is stated below as Theorem \ref{T:quant_CLT}. The proof of the quantile CLT is a straightforward exercise in the application of the Lindeberg-Feller theorem, and is given in the appendix.


To 
state
the quantile CLT, we first need some preliminaries. Let $X$ be an $\RR^d$-valued random variable. It will be convenient to denote vector components with function notation, $X=(X(1),\ldots,X(d))$. Let $\Phi_j(x) =P(X(j)\le x)$ and $G_{ij}(x,y)=P(X(i)\le x,X(j)\le y)$. Fix $\al\in(0,1)^d$ and assume there exists $q\in\RR^d$ such that $\Phi_j(q(j))=\al(j)$ for all $j$. Also assume that, for all $i$ and $j$, $\Phi_j'(q(j))$ exists and is strictly positive, and that $G_{ij}$ is continuous at $(q(i),q(j))$.

Let $\{X_n\}$ be a sequence of iid copies of $X$. If $\ka\in\{1,\ldots,n\}^d$, then $X_{\ka:n}\in\RR^d$ shall denote the componentwise order statistics of $X_1,\ldots,X_n$. That is, $X_{\ka:n}(j)$ is the $\ka(j)$-th order statistic of $(X_1(j),\ldots,X_n(j))$.

\begin{thm}\label{T:quant_CLT}
With the above notation and assumptions, define the matrix $\si\in\RR^{d\times d}$ by
  \begin{equation}\label{quant_CLT}
  \si_{ij} = \frac{G_{ij}(q(i),q(j)) - \al(i)\al(j)}
    {\Phi_i'(q(i))\Phi_j'(q(j))}.
  \end{equation}
If
$\ka(n)=(\ka(n,1),\ldots,\ka(n,d))\in\{1,\ldots,n\}^d$ satisfies $\ka(n)/n=\al+o(n^{-1/2})$, then $n^{1/2}(X_{\ka(n):n}-q)\To N$, where $N$ is multinormal with mean zero and covariance matrix $\si$.
\end{thm}

\begin{corollary}\label{C:quant_CLT}
There exists a centered Gaussian process $F$ with covariance function $\rho$, given by \eqref{rho_def}, such that $F_n\to F$ in the sense of finite-dimensional distributions.
\end{corollary}

\pf Given $0\le t_1<\cdots<t_d$, apply Theorem \ref{T:quant_CLT} with $X =(B(t_1),\ldots,B(t_d))$ and $\ka(n)=(j(n),j(n),\ldots,j(n))$. \qed

\begin{remark}
Note that one of the assumptions in the above theorem is that $\Phi_j'(q(j))$ exists and is strictly positive. In particular, $\Phi_j$ must be continuous at $q(j)$. That is, the theorem does not apply when the distribution of $X(j)$ has a jump at the quantile point $q(j)$. This, however, is never the case for the random variables considered in this paper, since for any initial distribution, the solution of the heat equation is continuous for any positive time.
\end{remark}

\begin{remark}
To give some limited intuitive sense to \eqref{quant_CLT}, let us consider what it implies for stochastic processes. Let $\{X_n(t): t\in T\}_{n=1}^\infty$ be a sequence of iid, measurable, real-valued copies of a stochastic process $X(t)$, parameterized by some general index set $T$. Let $\Phi(x,t) =P(X(t)\le x)$. Suppose that $\Phi(\cdot,t)$ is continuous and strictly increasing, and, for $\al\in(0,1)$, define $q(\al,t)$ so that $\Phi(q(\al,t),t)=\al$. Let $\ka:(0,1)\times\NN\to\NN$ satisfy $\ka(\al,n) \in\{1,\ldots,n\}$ and $\ka(\al,n)=\al n+o(n^{1/2})$ for each fixed $\al$. Let $X_{\ka:n}(t)$ denote the $\ka$-th order statistic of $X_1(t), \ldots,X_n(t)$, and define $Q_n(\al,t) =X_{\ka(\al,n):n} (t)$. Let  $F_n(\al,t)=n^{1/2} (Q_n(\al,t)-q(\al,t))$.

From Theorem \ref{T:quant_CLT}, we find that if we make the additional assumptions that $\pa_x\Phi(x,t)>0$ for all $x$ such that $0<\Phi(x,t)<1$, and that the functions $(x,y)\mapsto P(X(s)\le x,X(t)\le y)$ are continuous for each fixed $s$ and $t$, then the fdds of the two-parameter processes $F_n$ converge to those of a centered Gaussian process $F$ on $(0,1)\times T$ with covariance function
  \begin{equation}\label{cov_spec}
  E[F(\al,s)F(\be,t)] = \frac{P(X(s)\le q(\al,s), X(t)\le q(\be,t))
    - \al\be}
    {\pa_x\Phi(q(\al,s),s)\pa_x\Phi(q(\be,t),t)}.
  \end{equation}
In this paper, we are considering not only the special case in which the processes are Brownian motions, but also the special case in which $\al$ is fixed. If we fix $t$, however, then a familiar process emerges. Namely, for fixed $t$, the process $\al\mapsto F(\al,t)$ has covariance function
  \begin{equation}\label{bbcov}
  E[F(\al,t)F(\be,t)] = \frac{\al \wedge \be - \al\be}
    {\pa_x\Phi(q(\al,t),t)\pa_x\Phi(q(\be,t),t)}.
  \end{equation}
In other words, $\pa_x\Phi(q(\cdot,t),t)F(\cdot,t)$ is a Brownian bridge. This phenomenon stems from the relationship between the quantiles and the empirical processes. That is, let $\Phi^n(x,t) =\frac1n\sum_{j=1}^n1_{\{X_j(t)\le x\}}$. Then
  \[
  n^{1/2}(\Phi^n(Q_n(\al,t),t) - \Phi(q(\al,t),t))
    = n^{1/2}(n^{-1}\ka(\al,n) - \al)
    = o(1).
  \]
Hence, under suitable assumptions on $\Phi(\cdot,t)$ and its derivatives, we may use a Taylor polynomial for $\Phi(\cdot,t)$ about $q(\al,t)$, obtaining
  \begin{multline}\label{qerel}
  n^{1/2}(\Phi^n(Q_n(\al,t),t) - \Phi(Q_n(\al,t),t))
    = -n^{1/2}(\Phi(Q_n(\al,t),t) - \Phi(q(\al,t),t))
    + o(1)\\
  = -\pa_x\Phi(q(\al,t),t)F_n(\al,t) + n^{-1/2}|F_n(\al,t)|^2 O(1)
    + o(1).
  \end{multline}
For fixed $t$, it is well-known that $n^{1/2} (\Phi^n(\cdot,t) - \Phi(\cdot,t))$ converges to $B^\circ(\Phi(\cdot,t))$, where $B^\circ$ is a Brownian bridge (see Billingsley \cite{Bi}, for example). Hence, the left-hand side of \eqref{qerel} converges to $B^\circ(\Phi(q(\al,t),t)) = B^\circ (\al)$, which explains the covariance function \eqref{bbcov}.

The connection to the empirical processes given by \eqref{qerel} is still valid, of course, even when $t$ is allowed to vary. In this case, \eqref{qerel} and \eqref{cov_spec} show that the two-parameter fluctuation processes $n^{1/2}(\Phi^n-\Phi)\to V$ in the fdd sense, where $V$ is centered Gaussian with covariance
  \begin{equation}\label{Vcov}
  E[V(x,s)V(y,t)]
    = P(X(s)\le x,X(t)\le y) - P(X(s)\le x)P(X(t)\le y).
  \end{equation}
At this point, it is interesting to compare this with the results in Martin-L{\"o}f \cite{Mar}. When $X$ is a Markov process, \cite{Mar} considers the convergence, in the fdd sense, of fluctuations such as these, although for a slightly different model in which we have $N\sim\text{Poisson}(n)$ particles instead of $n$ particles. This subtle difference produces slightly different fdds in the limit. (For example, in the fixed $t$ case, the fluctuations of the empirical distributions of a Poisson number of variables converges to a Brownian motion rather than a Brownian bridge, which can be seen either from the heuristic demonstration below, or as a special case of Corollary 2 in \cite{Mar}.)

More specifically, let $\wt\Phi^n(x,t)=\frac1n\sum_{j=1}^N1_{\{X_j(t)\le x\}}$. Note that
  \[
  n^{1/2}(\wt\Phi^n - \Phi)
    = (N/n)^{1/2}N^{1/2}(\Phi^N - \Phi) + n^{1/2}(N/n - 1)\Phi.
  \]
Hence, $n^{1/2}(\wt\Phi^n-\Phi)\to\wt V$ in the fdd sense, where $\wt V = V + Z\Phi$, and $Z$ is a standard normal random variable, independent of $V$. Using \eqref{Vcov}, it follows that $\wt V$ is centered Gaussian with covariance
  \[
  E[\wt V(x,s)\wt V(y,t)] = P(X(s)\le x,X(t)\le y).
  \]
This, of course, is precisely the result we obtain from Corollary 2 in \cite{Mar}.
\end{remark}

\section{Properties of the limit process}\label{S:limit_props}

Before proceeding to the proof of tightness, we first establish some of the properties of the limit process. In what follows, $C$ shall denote a positive, finite constant, that may change value from line to line.

\begin{thm}\label{T:IV}
Recall $\rho$ given by \eqref{rho_def}. For each $T>0$, there exist positive constants $\de,C_1,C_2$ such that 
whenever $0<s<t\le T$ and $|t-s|<\de$, we have
  \begin{enumerate}[(i)]
  \item $C_1|t - s|^{-1/2} \le \pa_s\rho(s,t) \le C_2|t - s|^{-1/2}$,
  \item $-C_2|t - s|^{-1/2} \le \pa_t\rho(s,t) \le -C_1|t - s|^{-1/2}$, and
  \item $-C_2|t - s|^{-3/2} \le \pa^2_{st}\rho(s,t) \le -C_1|t - s|^{-3/2}$.
  \end{enumerate}
\end{thm}

\pf Define $\wt F(t)=u(q(t),t)F(t)$ and note that $\wt F$ is a centered Gaussian process with covariance function
  \begin{equation}\label{rho_tilde}
  \wt\rho(s,t) = P(B(s) \le q(s), B(t) \le q(t)) - \al^2.
  \end{equation}
We will first prove the theorem for $\wt\rho$ instead of $\rho$. For later purposes, it will be convenient at this time to introduce the auxiliary function
  \begin{equation}\label{zeta_def}
  \ze(s,t,w,z) = \int_{-\infty}^{q(s)+w}\int_{-\infty}^{q(t)+z}
    p(t-s,x,y)\pa_x^iu(x,s)\,dy\,dx,
  \end{equation}
where $i\ge0$ and $p(t,x,y)=(2\pi t)^{-1/2}e^{-(x-y)^2/2t}$.
In this proof, we use only the case $i=0$. Later, in Section \ref{S:param_est}, we will make use of the general case for arbitrary $i$. We now compute that
  \begin{equation}\label{newIV.1}
  \begin{split}
  \pa_t\ze &= q'(t)\int_{-\infty}^{q(s)+w}
    p(t-s,x,q(t)+z)\pa_x^i u(x,s)\,dx\\
  &\qquad + \int_{-\infty}^{q(s)+w}\int_{-\infty}^{q(t)+z}
    \pa_t p(t-s,x,y)\pa_x^i u(x,s)\,dy\,dx\\
  &= q'(t)\int_{-\infty}^{q(s)+w}
    p(t-s,x,q(t)+z)\pa_x^i u(x,s)\,dx\\
  &\qquad  + \frac12\int_{-\infty}^{q(t)+z}\int_{-\infty}^{q(s)+w}
    \pa_x^2 p(t-s,x,y)\pa_x^i u(x,s)\,dx\,dy.
  \end{split}
  \end{equation}
Using integration by parts,
  \begin{multline*}
  \int_{-\infty}^{q(s)+w}\pa_x^2 p(t-s,x,y)\pa_x^i u(x,s)\,dx\\
  = \pa_x p(t-s,q(s)+w,y)\pa_x^i u(q(s)+w,s)
    - p(t-s,q(s)+w,y)\pa_x^{i+1} u(q(s)+w,s)\\
    + \int_{-\infty}^{q(s)+w}p(t-s,x,y)\pa_x^{i+2} u(x,s)\,dx.
  \end{multline*}
Substituting this into \eqref{newIV.1} gives
  \begin{multline*}
  \pa_t\ze = q'(t)\int_{-\infty}^{q(s)+w}
    p(t-s,x,q(t)+z)\pa_x^i u(x,s)\,dx\\
  + \frac12\int_{-\infty}^{q(t)+z}
    \pa_x p(t-s,q(s)+w,y)\pa_x^i u(q(s)+w,s)\,dy\\
  - \frac12\int_{-\infty}^{q(t)+z}
    p(t-s,q(s)+w,y)\pa_x^{i+1} u(q(s)+w,s)\,dy\\
  + \frac12\int_{-\infty}^{q(t)+z}\int_{-\infty}^{q(s)+w}
    p(t-s,x,y)\pa_x^{i+2} u(x,s)\,dx\,dy.
  \end{multline*}
We now use the fact that $\pa_xp(t,x,y)=-\pa_yp(t,x,y)$ to rewrite this as
  \begin{multline}\label{newIV.2}
  \pa_t\ze = q'(t)\int_{-\infty}^{q(s)+w}
    p(t-s,x,q(t)+z)\pa_x^i u(x,s)\,dx\\
  - \frac12 p(t-s,q(s)+w,q(t)+z)\pa_x^i u(q(s)+w,s)\\
  - \frac12\pa_x^{i+1} u(q(s)+w,s)
    \int_{-\infty}^{q(t)+z}p(t-s,q(s)+w,y)\,dy\\
  + \frac12\int_{-\infty}^{q(t)+z}\int_{-\infty}^{q(s)+w}
    p(t-s,x,y)\pa_x^{i+2} u(x,s)\,dx\,dy.
  \end{multline}
Similarly, for the other partial derivative, we have
  \begin{multline*}
  \pa_s\ze = q'(s)\int_{-\infty}^{q(t)+z}
    p(t-s,q(s)+w,y)\pa_x^i u(q(s)+w,s)\,dy\\
    - \int_{-\infty}^{q(s)+w}\int_{-\infty}^{q(t)+z}
    \pa_t p(t-s,x,y)\pa_x^i u(x,s)\,dy\,dx\\
    + \int_{-\infty}^{q(s)+w}\int_{-\infty}^{q(t)+z}
    p(t-s,x,y)\pa_t\pa_x^i u(x,s)\,dy\,dx.
  \end{multline*}
By \eqref{newIV.1}, this becomes
  \begin{multline}\label{newIV.5}
  \pa_s\ze = q'(s)\int_{-\infty}^{q(t)+z}
    p(t-s,q(s)+w,y)\pa_x^i u(q(s)+w,s)\,dy\\
  - \pa_t\ze
    + q'(t)\int_{-\infty}^{q(s)+w} p(t-s,x,q(t)+z)\pa_x^i u(x,s)\,dx\\
  + \int_{-\infty}^{q(s)+w}\int_{-\infty}^{q(t)+z}
    p(t-s,x,y)\pa_t\pa_x^i u(x,s)\,dy\,dx.
  \end{multline}
Substituting \eqref{newIV.2} into the above, and using $\pa_tu=(1/2)\pa_x^2u$ gives
  \begin{multline}\label{newIV.3}
  \pa_s\ze = \frac12 p(t-s,q(s)+w,q(t)+z)\pa_x^i u(q(s)+w,s)\\
  + \left({q'(s)\pa_x^i u(q(s)+w,s)
    + \frac12\pa_x^{i+1} u(q(s)+w,s)}\right)
    \int_{-\infty}^{q(t)+z} p(t-s,q(s)+w,y)\,dy.
  \end{multline}
Now observe that, taking $i=0$, $\wt\rho(s,t) = \ze(s,t,0,0) - \al^2$. Hence, by \eqref{newIV.2} and \eqref{quant_ODE},
  \begin{multline*}
  \pa_t\wt\rho(s,t) = q'(t)\int_{-\infty}^{q(s)}
    p(t-s,x,q(t))u(x,s)\,dx
    - \frac12 p(t-s,q(s),q(t))u(q(s),s)\\
  + u(q(s),s)q'(s)\int_{-\infty}^{q(t)}p(t-s,q(s),y)\,dy
    + \frac12\int_{-\infty}^{q(t)}\int_{-\infty}^{q(s)}
    p(t-s,x,y)\pa_x^2 u(x,s)\,dx\,dy,
  \end{multline*}
and by \eqref{newIV.3} and \eqref{quant_ODE},
  \begin{equation}\label{IV.3}
  \pa_s\wt\rho(s,t) = \frac12 p(t-s,q(s),q(t))u(q(s),s).
  \end{equation}
Differentiating \eqref{IV.3} with respect to $t$ gives
  \begin{equation}\label{IV.4}
  \pa_{st}^2\wt\rho(s,t) = \frac12\pa_t p(t-s,q(s),q(t))u(q(s),s)
    + \frac12\pa_y p(t-s,q(s),q(t))u(q(s),s)q'(t).
  \end{equation}
Part (i) now follows easily from \eqref{IV.3}, using the fact that $u(q(s),s)$ is continuous and strictly positive on $[0,\infty)$, and $|q(t)-q(s)|\le C|t-s|$ for all $s,t\in[0,T]$.

Part (ii) will follow from (i), once we show that
  \begin{equation}\label{IV.9}
  |\pa_s\wt\rho(s,t) + \pa_t\wt\rho(s,t)| \le C.
  \end{equation}
By \eqref{newIV.5},
  \begin{multline}\label{IV.5}
  \pa_s\wt\rho(s,t) + \pa_t\wt\rho(s,t) = q'(s)\int_{-\infty}^{q(t)}
    p(t-s,q(s),y)u(q(s),s)\,dy\\
  + q'(t)\int_{-\infty}^{q(s)} p(t-s,x,q(t))u(x,s)\,dx
    + \int_{-\infty}^{q(s)}\int_{-\infty}^{q(t)}
    p(t-s,x,y)\pa_t u(x,s)\,dy\,dx.
  \end{multline}
Note that
  \begin{equation}\label{IV.6}
  \left|{q'(s)\int_{-\infty}^{q(t)}
    p(t-s,q(s),y)u(q(s),s)\,dy}\right|
    = |q'(s)|u(q(s),s)P^{q(s)}(B(t-s) \le q(t)) \le C.
  \end{equation}
Next, by the Markov property, $u(y,t) = \int_\RR p(t-s,x,y)u(x,s)\,dx$. Hence,
  \begin{equation}\label{IV.7}
  \left|{q'(t)\int_{-\infty}^{q(s)}p(t-s,x,q(t))u(x,s)\,dx}\right|
    \le |q'(t)|u(q(t),t) \le C.
  \end{equation}
Finally, note that $\pa_t u(x,t)=(1/2)\pa_x^2u(x,t)=(1/2)E^x[f''(B(t))]$. Therefore,
  \begin{multline*}
  \int_{-\infty}^{q(s)}\int_{-\infty}^{q(t)}
    p(t-s,x,y)\pa_t u(x,s)\,dy\,dx\\
  = \frac12\int_{-\infty}^{q(s)}\int_{-\infty}^{q(t)}\int_\RR
    p(t-s,x,y)p(s,x,z)f''(z)\,dz\,dy\,dx\\
  = \frac12\int_\RR P^z(B(s)\le q(s),B(t)\le q(t))f''(z)\,dz,
  \end{multline*}
which implies
  \begin{equation}\label{IV.8}
  \left|{\int_{-\infty}^{q(s)}\int_{-\infty}^{q(t)}
    p(t-s,x,y)\pa_t u(x,s)\,dy\,dx}\right|
    \le \frac12\|f''\|_1.
  \end{equation}
Applying \eqref{IV.6}, \eqref{IV.7}, and \eqref{IV.8} to \eqref{IV.5} establishes \eqref{IV.9}, and completes the proof of (ii).

Part (iii) will 
 follow
from \eqref{IV.4}, once we establish that for $|t-s|$ sufficiently small,
  \begin{align}
  -C_2|t - s|^{-3/2} \le \pa_t p(t-s,q(s),q(t))
    &\le -C_1|t - s|^{-3/2},\text{ and}\label{IV.10}\\
  |\pa_y p(t-s,q(s),q(t))| &\le C_2|t - s|^{-1/2}.\label{IV.11}
  \end{align}
First, note that
  \[
  \pa_t p(t,x,y) = -\frac1{2\sqrt{2\pi}}e^{-(x-y)^2/2t}
    \left({1 - \frac{(x-y)^2}t}\right)\frac1{t^{3/2}}.
  \]
Since $|q(t)-q(s)|\le C|t-s|$, it follows that for $|t-s|$ sufficiently small, $\pa_t p(t-s,q(s),q(t))<0$ with $|\pa_t p(t-s,q(s),q(t))|\le C|t-s|^{-3/2}$, which establishes the lower bound in \eqref{IV.10}. For the upper bound, note that $e^{-(q(t)-q(s))^2/2(t-s)}\ge e^{-C|t-s|}\ge e^{-CT}$. Hence, $|\pa_t p(t-s,q(s),q(t))| \ge C|t-s|^{-3/2}$ when $|t-s|$ is small enough.

Finally, we observe that
  \[
  |\pa_y p(t,x,y)| = \frac1{\sqrt{2\pi}}\frac{|x-y|}{t^{3/2}}
    e^{-(x-y)^2/2t} \le |x - y|t^{-3/2},
  \]
so that \eqref{IV.11} follows as above from $|q(t)-q(s)|\le C|t-s|$. This completes the proof of the theorem when $\rho$ is everywhere replaced by $\wt\rho$. We now prove the theorem as stated.

By \eqref{rho_def}, $\rho(s,t)=\th(s)\th(t)\wt\rho(s,t)$, where $\th(t) = (u(q(t),t))^{-1}$. Note that $\th\in C^1[0,\infty)$ and
  \begin{align*}
  \pa_s\rho(s,t) &= \th'(s)\th(t)\wt\rho(s,t)
    + \th(s)\th(t)\pa_s\wt\rho(s,t),\\
  \pa_t\rho(s,t) &= \th(s)\th'(t)\wt\rho(s,t)
    + \th(s)\th(t)\pa_t\wt\rho(s,t),\\
  \pa_{st}^2\rho(s,t) &= \th'(s)\th'(t)\wt\rho(s,t)
    + \th'(s)\th(t)\pa_t\wt\rho(s,t)
    + \th(s)\th'(t)\pa_s\wt\rho(s,t)
    + \th(s)\th(t)\pa_{st}^2\wt\rho(s,t),
  \end{align*}
from which (i), (ii), and (iii) follow. \qed

\begin{corollary}\label{C:var_order_0}
Fix $T>0$. There exist positive constants $C_1,C_2$ such that
  \begin{equation}\label{var_order_0}
  C_1|t - s|^{1/2} \le E|F(t) - F(s)|^2 \le C_2|t - s|^{1/2}
  \end{equation}
for all $s,t\in[0,T]$.
\end{corollary}

\pf Since $(s,t)\mapsto E|F(t)-F(s)|^2$ is continuous and strictly positive on $\{s\ne t\}$, it will suffice to show there exists $\de_0>0$ such that the theorem holds whenever $|t-s|\le\de_0$. Also, in the notation of the proof of Theorem \ref{T:IV}, since $F(t)=\th(t) \wt F(t)$, where $\th\in C^1 [0,\infty)$ is strictly positive, it will suffice to prove the theorem for $\wt F$. For this, observe that
  \begin{equation}\label{tilde_increment}
  E|\wt F(t) - \wt F(s)|^2
    = \wt\rho(t,t) + \wt\rho(s,s) - 2\wt\rho(s,t)
    = 2(\wt\rho(t,t) - \wt\rho(s,t)),
  \end{equation}
where $\wt\rho$ is given by \eqref{rho_tilde}. Hence,
  \[
  E|\wt F(t) - \wt F(s)|^2
    = 2\int_s^t \pa_s\wt\rho(u,t)\,du.
  \]
Applying Theorem \ref{T:IV}(i) and the fact that $\int_s^t |t-u|^{-1/2} \,du=2|t-s|^{1/2}$ finishes the proof. \qed

\begin{corollary}\label{C:Holder}
The process $F$ has a modification which is locally H\"older continuous on $[0,\infty)$ with exponent $\ga$ for any $\ga<1/4$.
\end{corollary}

\pf By the Kolmogorov-\v Centsov theorem (see, for example, \cite{KSh}, Thm 2.2.8), if, for each $T>0$, there exist positive constants $\be,r,C$ such that
  \[
  E|F(t) - F(s)|^\be \le C|t - s|^{1+r}
  \]
for all $s,t\in[0,T]$, then $F$ has a continuous modification which is locally H\"older continuous with exponent $\ga$ for all $\ga<r/\be$. Since $F$ is Gaussian, \eqref{var_order_0} implies $E|F(t)-F(s)|^p\le C_p|t-s|^{p/4}$. We may therefore take $\be=p$ and $r = p/4-1$. Letting $p\to\infty$ completes the proof. \qed

\begin{corollary}\label{C:IV.1}
For each $T>0$, there exist positive constants $\de,C_1,C_2$ such that
  \[
  -C_2|t - s|^{-1/2}\De t \le E[F(s)(F(t+\De t) - F(t))]
    \le -C_1|t + \De t - s|^{-1/2}\De t,
  \]
for all $0\le s<t\le T$ and $\De t>0$ such that $|t-s|<\de$ and $\De t<\de$.
\end{corollary}

\pf We write
  \[
  E[F(s)(F(t+\De t) - F(t))] = \rho(s,t+\De t) - \rho(s,t)
    = \int_t^{t+\De t} \pa_t\rho(s,u)\,du.
  \]
Applying Theorem \ref{T:IV}(ii) and the fact that
  \[
  |t + \De t - s|^{-1/2}\De t
    \le \int_t^{t+\De t} |u - s|^{-1/2}\,du
    \le |t - s|^{-1/2}\De t,
  \]
finishes the proof. \qed

\begin{corollary}\label{C:IV.2}
For each $T>0$, there exist positive constants $\de,C_1,C_2$ such that
  \begin{multline*}
  -C_2|t - s|^{-3/2}\De s\De t
    \le E[(F(s) - F(s-\De s))(F(t+\De t) - F(t))]\\
    \le -C_1|t + \De t - (s - \De s)|^{-3/2}\De s\De t,
  \end{multline*}
for all $0\le s<t\le T$, $\De s>0$, and $\De t>0$ such that $|t-s|<\de$, $\De s<\de$, and $\De t<\de$.
\end{corollary}

\pf We write
  \[
  E[(F(s) - F(s-\De s))(F(t+\De t) - F(t))]
    = \int_{s-\De s}^s\int_t^{t+\De t} \pa^2_{st}\rho(u,v)\,dv\,du.
  \]
Applying Theorem \ref{T:IV}(iii) and the fact that
  \[
  |t + \De t - (s - \De s)|^{-3/2}\De s\De t
    \le \int_{s-\De s}^s\int_t^{t+\De t} |v - u|^{-3/2}\,dv\,du
    \le |t - s|^{-3/2}\De s\De t,
  \]
finishes the proof. \qed

\begin{remark}\label{R:Ftilde}
Note that in Corollaries \ref{C:IV.1} and \ref{C:IV.2}, since $\pa_t \rho(s,t)$ and $\pa_{st}^2\rho(s,t)$ are continuous away from $\{s=t\}$, we have
  \begin{align*}
  |E[F(s)(F(t+\De t) - F(t))]|
    &\le C|t - s|^{-1/2}\De t,\text{ and}\\
  |E[(F(s) - F(s-\De s))(F(t+\De t) - F(t))]|
    &\le C|t - s|^{-3/2}\De s\De t,
  \end{align*}
for all $0\le s<t\le T$, $\De s>0$, and $\De t>0$. Also note that Theorem \ref{T:IV} and Corollaries \ref{C:var_order_0}-\ref{C:IV.2} all apply to $\wt F(t)=u(q(t),t)F(t)$ as well.
\end{remark}

As an application of these estimates on the local covariance structure of the increments of $F$, we now show that $F$ has a finite, nonzero, deterministic quartic variation.

\begin{thm}\label{T:quart_var}
If
  \[
  V_\Pi(t) = \sum_{0<t_j\le t} |F(t_j) - F(t_{j-1})|^4,
  \]
where $\Pi=\{0=t_0<t_1<t_2<\cdots\}$ is a partition of $[0,\infty)$ with $t_j\uparrow\infty$, then for all $T>0$,
  \[
  \lim_{|\Pi|\to0} E\bigg[\sup_{0\le t\le T}
    \bigg|V_\Pi(t)
    - \frac6\pi\int_0^t|u(q(s),s)|^{-2}\,ds\bigg|^2\bigg]
    = 0,
  \]
where $|\Pi|=\sup_j|t_j-t_{j-1}|$.
\end{thm}

\pf We again adopt the notation of the proof of Theorem \ref{T:IV}. Let $\th(t)=(u(q(t),t))^{-1}$ and recall that $\th\in C^1[0,\infty)$ is strictly positive, and that $F(t)=\th(t)\wt F(t)$. Since $V_\Pi$ is monotone, it suffices by Dini's theorem to show that
  \[
  V_\Pi(t) \to \frac6\pi\int_0^t\th(s)^2\,ds,
  \]
in $L^2$ for each fixed $t>0$. For any stochastic process $X$, we adopt the notation $\De X_j=X(t_j)-X(t_{j-1})$. We will also write $\De t_j=t_j-t_{j-1}$. Hence,
  \[
  \De F_j = \th(t_j)\De\wt F_j + \wt F(t_{j-1})\De\th_j,
  \]
which implies $\De F_j^4 = \th(t_j)^4\De\wt F_j^4 + R_j$, where
  \[
  |R_j| \le C\sum_{i=1}^4 |\wt F(t_{j-1})|^i\De t_j^i\De\wt F_j^{4-i}.
  \]
Thus,
  \[
  E\bigg|\sum_{0<t_j\le t} R_j\bigg|^2 \le C\sum_{i=1}^4
    E\bigg|\sum_{0<t_j\le t}
    |\wt F(t_{j-1})|^i\De t_j^i\De\wt F_j^{4-i}\bigg|^2.
  \]
Note that Corollary \ref{C:var_order_0} also holds for $\wt F$. Thus, using H\"older's inequality, we have
  \begin{align*}
  E\bigg|\sum_{0<t_j\le t} R_j\bigg|^2
    &\le C\sum_{i=1}^4 E\bigg[\bigg(\sum_{0<t_j\le t}\De t_j^2\bigg)
    \bigg(\sum_{0<t_j\le t}
    |\wt F(t_{j-1})|^{2i}\De t_j^{2(i-1)}\De\wt F_j^{2(4-i)}
    \bigg)\bigg]\\
  &= C\sum_{i=1}^4 \bigg(\sum_{0<t_j\le t}\De t_j^2\bigg)
    \bigg(\sum_{0<t_j\le t}
    \De t_j^{2(i-1)}E[|\wt F(t_{j-1})|^{2i}\De\wt F_j^{2(4-i)}]
    \bigg)\\
  &\le C\sum_{i=1}^4 \bigg(\sum_{0<t_j\le t}\De t_j^2\bigg)
    \bigg(\sum_{0<t_j\le t}
    \De t_j^{2(i-1)}(E|\wt F(t_{j-1})|^{4i})^{1/2}
    (E\De\wt F_j^{4(4-i)})^{1/2}
    \bigg)\\
  &\le C\sum_{i=1}^4 \bigg(\sum_{0<t_j\le t}\De t_j^2\bigg)
    \bigg(\sum_{0<t_j\le t}
    \De t_j^{2i-2}(E\De\wt F_j^{4(4-i)})^{1/2}
    \bigg).
  \end{align*}
As in the proof of Corollary \ref{C:Holder}, since $F$ is Gaussian, \eqref{var_order_0} implies $E\De\wt F_j^{4(4-i)} \le C_i\De t^{4-i}$. Thus,
  \[
  E\bigg|\sum_{0<t_j\le t} R_j\bigg|^2
    \le C\sum_{i=1}^4 \bigg(\sum_{0<t_j\le t} \De t_j^2\bigg)
    \bigg(\sum_{0<t_j\le t} \De t_j^{3i/2}\bigg) \to 0,
  \]
as $|\Pi|\to0$. It therefore suffices to show that
  \[
  \sum_{0<t_j\le t} \th(t_j)^4\De\wt F_j^4
    \to \frac6\pi\int_0^t\th(s)^2\,ds,
  \]
in $L^2$ as $|\Pi|\to0$.

By \eqref{tilde_increment},
  \begin{equation}\label{quart_var.1}
  E\De\wt F_j^2 = 2(\wt\rho(t_j,t_j) - \wt\rho(t_{j-1},t_j))
    = 2\int_0^{\De t_j}\pa_s\wt\rho(t_j-\ep,t_j)\,d\ep.
  \end{equation}
By \eqref{IV.3},
  \begin{multline*}
  \pa_s\wt\rho(t_j-\ep,t_j) = \frac12 p(\ep,q(t_j-\ep),q(t_j))
    u(q(t_j-\ep),t_j-\ep)\\
  = \frac1{2\sqrt{2\pi}}\,\ep^{-1/2}e^{-(q(t_j)-q(t_j-\ep))^2/2\ep}
    u(q(t_j-\ep),t_j-\ep).
  \end{multline*}
Using $|1-e^{-x}|\le x$ for all $x\ge0$, and $|q(t)-q(s)|\le C|t-s|$ for all $0\le s\le t\le T$, and $u(q(\cdot),\cdot)\in C^1[0,\infty)$, this gives
  \[
  \left|\pa_s\wt\rho(t_j-\ep,t_j) - \frac1{2\sqrt{2\pi}}\,\ep^{-1/2}
    u(q(t_j),t_j)\right| \le C\ep^{1/2}.
  \]
Substituting this into \eqref{quart_var.1} gives
  \[
  \bigg|E\De\wt F_j^2
    - \sqrt{\frac2\pi}\De t_j^{1/2}u(q(t_j),t_j)\bigg| \le C\De t_j^{3/2}.
  \]
Since $E\De\wt F_j^4=3(E\De\wt F_j^2)^2$ and $\th(t)=(u(q(t),t))^{-1}$, this implies
  \begin{multline*}
  \bigg|\th(t_j)^4 E\De\wt F_j^4
    - \frac6\pi\De t_j \th(t_j)^2\bigg|\\
    = 3\th(t_j)^4\bigg|E\De\wt F_j^2
    - \sqrt{\frac2\pi}\De t_j^{1/2} \th(t_j)^{-1}\bigg|
    \bigg|E\De\wt F_j^2
    + \sqrt{\frac2\pi}\De t_j^{1/2} \th(t_j)^{-1}\bigg|
    \le C\De t_j^2.
  \end{multline*}
Note that
  \[
  \sum_{0<t_j\le t} \frac6\pi\De t_j \th(t_j)^2
    \to \frac6\pi\int_0^t \th(s)^2\,ds.
  \]
It will therefore suffice to show that
  \[
  E\bigg|\sum_{0<t_j\le t} (\th(t_j)^4\De\wt F_j^4
    - \th(t_j)^4 E\De\wt F_j^4)\bigg|^2 \to 0,
  \]
as $|\Pi|\to0$.

For this, we write
  \[
  E\bigg|\sum_{0<t_j\le t} (\th(t_j)^4\De\wt F_j^4
    - \th(t_j)^4 E\De\wt F_j^4)\bigg|^2
    = \sum_{i,j} \th(t_i)^4\th(t_j)^4(E[\De\wt F_i^4\De\wt F_j^4]
    - (E\De\wt F_i^4)(E\De\wt F_j^4)).
  \]
If $X$ and $Y$ are jointly normal with mean zero and variances $\si_X^2$ and $\si_Y^2$, then
  \[
  E[X^4Y^4] = \si_X^4\si_Y^4(9 + 72\rho_{X,Y}^2 + 24\rho_{X,Y}^4),
  \]
where $\rho_{X,Y} = (\si_X\si_Y)^{-1}E[XY]$. Hence,
  \[
  |E[X^4Y^4] - (EX^4)(EY^4)| \le C\si_X^2\si_Y^2|E[XY]|^2.
  \]
Therefore,
  \[
  E\bigg|\sum_{0<t_j\le t} (\th(t_j)^4\De\wt F_j^4
    - \th(t_j)^4 E\De\wt F_j^4)\bigg|^2
    \le C\sum_{i,j}\De t_i^{1/2}\De t_j^{1/2}|E\De\wt F_i\De\wt F_j|^2.
  \]
By H\"older's inequality, $|E\De\wt F_i\De\wt F_j|^2\le C\De t_i^{1/2} \De t_j^{1/2}$. Thus, it will suffice to show that
  \[
  \sum_{t_j-t_i>2|\Pi|}
    \De t_i^{1/2}\De t_j^{1/2}|E\De\wt F_i\De\wt F_j|^2
    \to 0,
  \]
as $|\Pi|\to0$.

If $t_j-t_i>2|\Pi|$, then $t_{j-1}>t_i$. Hence, by Corollary \ref{C:IV.2} and Remark \ref{R:Ftilde},
  \[
  |E\De\wt F_i\De\wt F_j|^2
    \le C\De t_i^2\De t_j^2|t_{j-1} - t_i|^{-3},
  \]
where $C$ depends only on $T$. Since $|t_{j-1}-t_i|>|\Pi|\ge\De t_k$ for all $k$, this implies
  \[
  |E\De\wt F_i\De\wt F_j|^2
    \le C\De t_i^{1/2}\De t_j^{5/4}|t_{j-1} - t_i|^{-3/4}
    \le C|\Pi|^{3/4}\De t_i^{1/2}\De t_j^{1/2}|t_{j-1} - t_i|^{-3/4}.
  \]
Hence,
  \[
  \sum_{t_j-t_i>2|\Pi|}
    \De t_i^{1/2}\De t_j^{1/2}|E\De\wt F_i\De\wt F_j|^2
    \le C|\Pi|^{3/4}\sum_{t_j-t_i>2|\Pi|}
    \De t_i\De t_j|t_{j-1} - t_i|^{-3/4} \to 0,
  \]
since $\int_0^t\int_0^t|u-v|^{-3/4}\,du\,dv<\infty$. \qed

\bigskip

Although the local properties of $F$ are very similar to $B^{1/4}$, the global properties need not be.
One simple observation along these lines is that $E|B^H(t)|^2=t^{2H}$, whereas $E|F(t)|^2\ge Ct$, for some constant $C$ that depends on the initial distribution $f$. Indeed, by \eqref{rho_def}, we have $E|F(t)|^2 =(\al-\al^2)|u(q(t),t)|^{-2}$, and
  \[
  u(x,t) = \frac1{\sqrt{2\pi t}}\int f(y)e^{-(x-y)^2/2t}\,dy
    \le \frac1{\sqrt{2\pi t}}\|f\|_1.
  \]
In other words, for $t$ large, $F(t)$ has a variance which is comparable to Brownian motion, rather than fBm.

We will illustrate some other global properties of $F$
 with a particular example of an initial distribution. As noted in Section \ref{S:intro}, if $B(0)\sim N(0,1)$  and $j(n) = \flr{{(n+1)}/2}$, so that $Q_n$ is the median, $\al=1/2$, and $q=0$, then the result in \cite{Sw07.2} implies that $F_n(\cdot)\To X(\cdot+1)$ in $C[0, \infty)$, where the covariance function of $X$ is given by \eqref{med_covar}. We will illustrate the global properties of $F$ in this case by deriving some of the global properties of $X$.

First, note that fBm has a self-similarity scaling property. Namely, $c^{-H}B^H(ct)$ and $B^H(t)$ have the same law, as processes. Based on this, we might expect $X$ to have a similar scaling property with exponent $1/4$. However, $X$ obeys the Brownian scaling law with exponent $1/2$. This can be seen immediately by direct inspection of \eqref{med_covar}.

Second, note that fBm is a long-memory process. In particular,
  \[
  r_H(n) = E[B^H(1)(B^H(n+1) - B^H(n))]
  \]
decays only polynomially. It is well-known that $r_H(n)\sim2H(2H-1)n^{2H-2}$, where $a_n\sim b_n$ means $a_n/b_n\to1$ as $n\to\infty$. Based on this, we might expect
  \begin{equation}\label{long_mem}
  r(n) = E[X(1)(X(n+1) - X(n))]
  \end{equation}
to satisfy $r(n)\sim Cn^{-3/2}$ for some constant $C$. However, this is not the case.

\begin{prop}\label{P:incr_decay}
If $r(n)$ is given by \eqref{long_mem}, then $r(n)\sim -(1/6)n^{-2}$.
\end{prop}

\pf The proposition follows easily by observing that
  \[
  r(n) = \sqrt{n+1}\,\sin^{-1}\left(\frac1{\sqrt{n+1}}\right)
    - \sqrt n\,\sin^{-1}\left(\frac1{\sqrt n}\right),
  \]
and using the Taylor expansion $\sin^{-1}x = x + x^3/6 + 3x^5/40 + O(x^7)$. \qed

\bigskip

We see then that $X$ is a process which behaves locally like fBm with $H=1/4$, but whose increments are more weakly correlated than those of fBm. Another example of such a process is sub-fBm (see Bojdecki et al \cite{BGT}, for example). For sub-fBm with $H=1/4$, the increments decay at the rate $n^{-5/2}$. The decay rate of $n^{-2}$ in Proposition \ref{P:incr_decay} has only been established here for the case when the initial distribution of the particles is a standard Gaussian, and an investigation of the relationship between the initial distribution and this exponent is beyond the scope of the current paper. However, we may presently note that $F$ is not a sub-fBm for any initial distribution $f$, since the variance of a sub-fBm is always proportional to $t^{2H}$, whereas $E|F(t)|^2\ge Ct$, as we observed above.

\section{Outline of proof of tightness}

Our primary tool for proving tightness is the following version of the Kolmogorov-Prohorov tightness criterion.

\begin{thm}\label{T:moment}
If $\{Z_n\}$ is a sequence of continuous stochastic processes such that
  \begin{enumerate}[(i)]
  \item $\sup_{n\ge n_0} P(|Z_n(0)| > \la) \to 0$ as $\la\to\infty$, and
  \item $\sup_{n\ge n_0} P(|Z_n(t) - Z_n(s)| \ge \ep) \le C_T\ep^{-\al}|t - s|^{1+\beta}$ whenever $0 < \ep < 1$, $T > 0$, and $0 \le s<t \le T$,
  \end{enumerate}
for some positive constants $n_0$, $\al$, $\beta$, and $C_T$ (depending on $T$), then $\{Z_n\}$ is relatively compact in $C[0,\infty)$.
\end{thm}

When we apply the above theorem to our processes $\{F_n\}$, Condition (i) will automatically be satisfied due to the already established convergence of the fdds.
Our
 main theorem will
therefore
 be proved once we establish that $\{F_n\}$ satisfies Condition (ii) of Theorem \ref{T:moment}. The remainder of the paper will be devoted to proving this. For this, we will need to study the distribution of the difference of two quantiles. Unfortunately, the difference of the quantiles is not the quantile of the differences. The stark non-linearity of the quantile mapping creates substantial technical difficulties. In the next section, we begin by using conditioning to connect the probability being estimated to a certain random walk. With an eye to future applications, we construct this connection in abstract generality. In Section \ref{S:param_est}, we will return to the specific sequence $\{F_n\}$.

Before moving on to the details of the proof of tightness, we pause here to present a general overview of the key steps to be taken in the remainder of the paper. Let $\ol B_j = B_j - q$, $\ol B = \ol B_1$, and $\ol Q_n = Q_n - q = \ol B_{j(n):n}$. To verify Condition (ii) of Theorem \ref{T:moment}, we will aim to show that
  \begin{equation}\label{heur_star2}
  P(|F_n(t) - F_n(s)| \ge \ep) \le C_p\ep^{-p}|t - s|^{p/4},
  \end{equation}
for any $p>2$. We begin by writing
  \[
  P(|F_n(t) - F_n(s)| \ge \ep)
    = P(|\ol Q_n(t) - \ol Q_n(s)| \ge n^{-1/2}\ep).
  \]
Heuristically, the order of magnitude of $|\ol Q_n(t) - \ol Q_n(s)|$ should be less than that of $|\ol B(t) - \ol B(s)|$, which is $|t - s|^{1/2}$. Consequently, it is easy to estimate this probability when $n^{-1/2}\ep\ge|t-s|^{1/2-\De}$ for some $\De>0$. (This is what we call the ``large gap regime".) This routine estimate is carried out in the appendix, and \eqref{heur_star2} is easily verified in the large gap regime.

We next turn our attention to the ``small gap regime", when $n^{-1/2}\ep\le|t-s|^{1/2+\De}$ for some $\De>0$. We begin with the fact that
  \begin{equation}\label{heur_star}
  P(\ol Q_n(t) - \ol Q_n(s) \le -n^{-1/2}\ep)
    = E[P(\ol Q_n(t) - \ol Q_n(s) \le -n^{-1/2}\ep
    \mid \ol Q_n(s))].
  \end{equation}
(The estimates for $P(\ol Q_n(t) - \ol Q_n(s) \ge n^{-1/2}\ep)$ are nearly identical.) To deal with the right-hand side of \eqref{heur_star}, we apply the results from Section \ref{S:RWrep}, wherein we study, in generality, the conditional distribution of the difference of two quantiles. There we find that this conditional probability can be well-approximated by certain probabilities connected to sums of iid random variables, and that sufficient estimates on these probabilities can be derived. When we apply Theorem \ref{T:RWrep} to the problem at hand, we find that if
  \[
  \begin{split}
  q_1(x,y) &= P(B(t) > q(t) + x + y \mid B(s) < q(s) + x),\\
  q_2(x,y) &= P(B(t) < q(t) + x + y \mid B(s) > q(s) + x),
  \end{split}
  \]
and
  \[
  \ph_{j:n}^\le(x,y) = P\bigg(\sum_{i=1}^{j-1} 1_{\{U_i\le q_1\}}
    \le \sum_{i=j+1}^n 1_{\{U_i\le q_2\}}\bigg),
  \]
where $\{U_i\}$ are iid and uniform on $(0,1)$, then
  \[
  P(\ol Q_n(t) - \ol Q_n(s) \le -n^{-1/2}\ep \mid \ol Q_n(s))
    \le \ph_{j(n):n}^\le(\ol Q_n(s),-n^{-1/2}\ep).
  \]
By \eqref{heur_star} and the fact that $\ph_{j:n}^\le(x,y)\le1$, this gives, for any $K>0$,
  \[
  P(\ol Q_n(t) - \ol Q_n(s) \le -n^{-1/2}\ep)
    \le \sup_{|x|\le K}\ph_{j(n):n}^\le(x,-n^{-1/2}\ep)
    + P(|\ol Q_n(s)|>K).
  \]
Since $F_n$ has well-behaved tail probabilities (see Proposition \ref{P:quant_tails} and Remark \ref{R:quant_tails}), we may choose $K$ so that $P(|\ol Q_n(s)|>K)\le C_p\ep^{-p}|t - s|^{p/4}$. (In fact, the correct choice is $K=n^{-1/2}\ep|t-s|^{-1/4}$.) The problem thus reduces to showing that
  \begin{equation}\label{heur_star3}
  \sup_{|x|\le K}\ph_{j(n):n}^\le(x,-n^{-1/2}\ep)
    \le C_p\ep^{-p}|t - s|^{p/4}.
  \end{equation}
For this, we again appeal to the general results in Section \ref{S:RWrep}. This time applying Theorem \ref{T:RWest}, we find that
  \[
  \ph_{j(n):n}^\le(x,-n^{-1/2}\ep)
    \le C\left({\frac{\si(x,-n^{-1/2}\ep)}
    {n\mu(x,-n^{-1/2}\ep)^2}
    }\right)^{p/2},
  \]
where
  \begin{align*}
  \si(x,y) &= \al q_1(x,y) + (1 - \al)q_2(x,y),\\
  \mu(x,y) &= \al q_1(x,y) - (1 - \al)q_2(x,y).
  \end{align*}
To show that this leads to \eqref{heur_star3}, we use Taylor expansions for $\si$ and $\mu$ which are developed in Proposition \ref{P:Psi_Taylor} and Corollary \ref{C:Psi_Taylor}, and which rely on some of the computations in the proof of Theorem \ref{T:IV}. When all of this analysis is finally complete, we will have established \eqref{heur_star3}, and therefore \eqref{heur_star2}, in the small gap regime.

Lastly, we must deal with the ``medium gap regime", in which $|t-s|^{1/2+\De} \le n^{-1/2}\ep\le |t-s|^{1/2-\De}$ for some $\De>0$. This case is treated by making minor modifications to the proof for the small gap regime. In making these modifications, however, the result is weakened, and we are only able to prove that, for parameter values in the medium gap regime,
  \begin{equation}\label{heur_star4}
  P(|F_n(t) - F_n(s)| > \ep)
    \le C(\ep^{-1}|t - s|^{1/4-2\De})^p.
  \end{equation}
Although we conjecture that the sharper bound \eqref{heur_star2} does in fact hold in all regimes, the weaker bound \eqref{heur_star4} which we are able to prove is still sufficient to establish tightness.

\section{A random walk representation and estimate}\label{S:RWrep}

Let $\{U_n\}_{n=1}^\infty$ be an iid sequence of random variables, uniformly distributed on $(0,1)$, and
for $r_1,r_2\in(0,1)$,
define
  \begin{align*}
  \psi_{j:n}^\le(r_1,r_2) = P\bigg(\sum_{i=1}^{j-1} 1_{\{U_i\le r_1\}}
    \le \sum_{i=j+1}^n 1_{\{U_i\le r_2\}}\bigg),\\
  \psi_{j:n}^<(r_1,r_2) = P\bigg(\sum_{i=1}^{j-1} 1_{\{U_i\le r_1\}}
    < \sum_{i=j+1}^n 1_{\{U_i\le r_2\}}\bigg).
  \end{align*}
Also, let $\psi_{j:n}^>=1-\psi_{j:n}^\le$ and $\psi_{j:n}^\ge=1-\psi_{j:n}^<$.

Let $X$ and $Y$ be real-valued random variables such that $(x,y) \mapsto P(X\le x,Y\le y)$ is continuous. Define
  \[
  \begin{split}
  q_1(x,y) &= P(Y > x + y \mid X < x),\\
  q_2(x,y) &= P(Y < x + y \mid X > x),
  \end{split}
  \]
and let $\ph_{j:n}^\star(x,y)=\psi_{j:n}^\star(q_1(x,y),q_2(x,y))$, where $\star$ may be any one of the symbols $\le$, $<$, $>$, or $\ge$. Note that $\ph_{j:n}^>(x,y) = \psi^<_{(n-j+1):n}(q_2,q_1)$ and $\ph^\ge_{j:n} (x,y) = \psi^\le_{(n-j+1):n}(q_2,q_1)$.

\begin{thm}\label{T:RWrep}
If $\{(X_n,Y_n)\}$ is a sequence of iid copies of $(X,Y)$, then for all $y\in\RR$,
  \begin{equation}\label{RWrep.1}
  \ph_{j:n}^<(X_{j:n},y) \le P(Y_{j:n} - X_{j:n} < y \mid X_{j:n})
    \le \ph_{j:n}^\le(X_{j:n},y),
  \end{equation}
almost surely. Consequently,
  \begin{equation}\label{RWrep.2}
  \ph_{j:n}^>(X_{j:n},y) \le P(Y_{j:n} - X_{j:n} > y \mid X_{j:n})
    \le \ph^\ge_{j:n}(X_{j:n},y),
  \end{equation}
almost surely.
\end{thm}

Theorem 5.1 establishes a connection between the conditional distribution of the difference of two quantiles and probabilities connected to sums of iid random variables. To give some intuitive sense to this connection, let us consider the following heuristic derivation of \eqref{RWrep.1}.

We are interested in estimating $P(Y_{j:n} - X_{j:n} < y \mid X_{j:n} = x) = P(Y_{j:n} < x + y \mid X_{j:n} = x)$. Let us consider $X_1,\ldots, X_n$ as representing the locations on the real line of some particles at time $s$, and $Y_1,\ldots,Y_n$ as representing the locations of those same particles at some later time $t>s$. We are given that the $j$-th leftmost particle at time $s$ is located at position $x$. Conditioned on this information, we know the following. At time $s$, there is one particle located at $x$, there are $j-1$ iid particles located in $(-\infty,x)$, and there are $n-j$ iid particles located in $(x,\infty)$. The event $\{Y_{j:n} < x + y\}$ will occur if and only if at least $j$ particles end up in $(-\infty,x+y)$ at time $t$.

Each particle which is in $(-\infty,x)$ at time $s$ has probability
  \[
  p_1 = 1 - q_1 = P(Y < x + y \mid X < x)
  \]
of ending up in the target interval $(-\infty,x+y)$ at time $t$. Therefore, we may represent the number of particles from $(-\infty, x)$ which end up in the target interval by $\sum_{i=1}^{j-1} 1_{\{ U_i\le p_1\}}$. Similarly, the number of particles from $(x,\infty)$ which end up in the target interval is represented by $\sum_{i=j+1}^n 1_{\{U_i\le q_2\}}$.

We are therefore interested in computing the probability that
  \begin{align*}
  j &\le \sum_{i=1}^{j-1} 1_{\{ U_i\le p_1\}}
    + \sum_{i=j+1}^n 1_{\{U_i\le q_2\}}\\
  &\overset{d}{=} \sum_{i=1}^{j-1} (1 - 1_{\{ U_i\le q_1\}})
    + \sum_{i=j+1}^n 1_{\{U_i\le q_2\}}\\
  &= j - 1 - \sum_{i=1}^{j-1} 1_{\{ U_i\le q_1\}}
    + \sum_{i=j+1}^n 1_{\{U_i\le q_2\}},
  \end{align*}
which happens if and only if $\sum_{i=1}^{j-1} 1_{\{ U_i\le q_1\}} < \sum_{i=j+1}^n 1_{\{U_i\le q_2\}}$. This probability is exactly $\psi_{j:n}^<(q_1,q_2)=\ph_{j:n}^<(x,y)$. In fact, the true probability is even larger than this, since the particle at $x$ may itself end up in the target interval $(-\infty,x+y)$ at time $t$. If that happens, then we only need $\sum_{i=1}^{j-1}1_{\{U_i\le q_1\}}\le\sum_{i=j+1}^n 1_{\{U_i\le q_2\}}$, and the probability of this is $\ph_{j:n}^\le(x,y)$. Hence, the conditional probability of interest is sandwiched between $\ph_{j:n}^<$ and $\ph_{j:n}^\le$.

\bigskip

\noindent{\bf Proof of Theorem \ref{T:RWrep}.} By taking complements, \eqref{RWrep.1} and \eqref{RWrep.2} are equivalent. We will prove \eqref{RWrep.2}. For this, we will first establish that for any $a,b,y\in\RR$ with $a<b$,
  \begin{multline}\label{RWrep.3}
  E[\ph_{j:n}^>(X_{j:n},y)1_{\{X_{j:n} \in (a,b)\}}]
    \le P(Y_{j:n} - X_{j:n} > y, X_{j:n} \in (a,b))\\
    \le E[\ph_{j:n}^\ge(X_{j:n},y)1_{\{X_{j:n} \in (a,b)\}}].
  \end{multline}
Let us begin with the upper bound in \eqref{RWrep.3}.

First, by taking $a_i\downarrow a$ and $b_i\uparrow b$ if necessary, we may assume that $P(X<a) >0$ and $P(X>b)>0$, so that $q_i$, $\ph_{j:n}^\ge$, and $\Phi$, where $\Phi(x)=P(X<x)$, are all well-defined, continuous, and bounded on $\{(x,y)\in[a,b]\times\RR\}$. Now fix $\ep>0$ and $x\in[a,b-\ep]$ and let $N=\#\{i:X_i\in(x,x+\ep)\}$. Then
  \begin{multline*}
  P(Y_{j:n} - X_{j:n} > y, X_{j:n} \in (x,x+\ep)) \le P(N \ge 2)\\
    + P(Y_{j:n} > x + y, \#\{i: X_i < x\} = j - 1,
    \#\{i: X_i > x + \ep\} = n - j, N = 1).
  \end{multline*}
Thus,
  \begin{multline}\label{RWrep.4}
  P(Y_{j:n} - X_{j:n} > y, X_{j:n} \in (x,x+\ep))\\
    \le \binom{n}{2}P(X \in (x,x+\ep))^2
    +  j\binom{n}{j}P(Y_{j:n} > x + y, A),
  \end{multline}
where
  \[
  A = \{X_i < x\text{ for all }i<j\} \cap
    \{X_i > x + \ep\text{ for all }i>j\} \cap \{X_j \in (x,x+\ep)\}.
  \]
Note that
  \begin{align*}
  P(Y_{j:n} > x + y,\, &A)\\
  &= P(\#\{i: Y_i < x + y\} < j, A)\\
  &\le P(\#\{i \ne j: Y_i < x + y\} < j, A)\\
  &= \sum_{\ell=0}^{j-1}\sum_{m=\ell}^{j-1}
    P(\#\{i < j: Y_i < x + y\} = \ell,
    \#\{i > j: Y_i < x + y\} = m - \ell, A).
  \end{align*}
By independence, this gives
  \begin{multline*}
  P(Y_{j:n} > x + y, A)\\
  \le \sum_{\ell=0}^{j-1}\sum_{m=\ell}^{j-1}
    \binom{j-1}{\ell}P(Y < x + y, X < x)^\ell
    P(Y > x + y, X < x)^{j-1-\ell}\\
  \cdot\binom{n-j}{m-\ell}
    P(Y < x + y, X > x + \ep)^{m-\ell}
    P(Y > x + y, X > x + \ep)^{n-j-m+\ell}\\
  \cdot P(X\in(x,x+\ep)).
  \end{multline*}
If we define $\wt q_i(x,y)=q_i(x+\ep,y-\ep)$, $p_i=1-q_i$, $\wt p_i = 1 - \wt q_i$, and $\ol\Phi=1-\Phi$, then this becomes
  \begin{multline}\label{RWrep.5}
  P(Y_{j:n} > x + y, A) \le \sum_{\ell=0}^{j-1}\sum_{m=\ell}^{j-1}
    \binom{j-1}{\ell}\binom{n-j}{m-\ell}
    p_1^\ell q_1^{j-1-\ell}\wt p_2^{n-j-m+\ell}\wt q_2^{m-\ell}\\
    \cdot\Phi(x)^{j-1}\ol\Phi(x+\ep)^{n-j}P(X\in(x,x+\ep)).
  \end{multline}
Compare this with
  \begin{align}
  \ph_{j:n}^\ge(x,y)
    &= P\bigg(\sum_{i=1}^{j-1} 1_{\{U_i\le q_1\}}
    \ge \sum_{i=j+1}^n 1_{\{U_i\le q_2\}}\bigg)\notag\\
  &= \sum_{\ell=0}^{j-1}\sum_{m=\ell}^{j-1}
    P\bigg(\sum_{i=1}^{j-1} 1_{\{U_i\le q_1\}} = j - 1 - \ell,
    \sum_{i=j+1}^n 1_{\{U_i\le q_2\}} = m - \ell\bigg)\notag\\
  &= \sum_{\ell=0}^{j-1}\sum_{m=\ell}^{j-1}
    \binom{j-1}{\ell}\binom{n-j}{m-\ell}
    p_1^\ell q_1^{j-1-\ell}p_2^{n-j-m+\ell}q_2^{m-\ell}.\label{RWrep.6}
  \end{align}
Finally, partition $(a,b)$ into subintervals of size $\ep$ and apply  \eqref{RWrep.4} and \eqref{RWrep.5}. Let $\ep\to0$ and apply dominated convergence. By \eqref{RWrep.6}, this gives
  \[
  P(Y_{j:n} - X_{j:n} > y, X_{j:n} \in (a,b))
    \le j\binom{n}{j}\int_a^b \ph_{j:n}^\ge(x,y)
    \Phi(x)^{j-1}\ol\Phi(x)^{n-j}\,d\Phi(x).
  \]
By Lemma \ref{L:RWrep}, this completes the proof of the upper bound in \eqref{RWrep.3}.

The lower bound in \eqref{RWrep.3} is proved similarly. As before,
  \begin{multline*}
  P(Y_{j:n} - X_{j:n} > y, X_{j:n} \in (x,x+\ep))\\
    \ge P(Y_{j:n} > x + \ep + y, \#\{i: X_i < x\} = j - 1,
    \#\{i: X_i > x + \ep\} = n - j, N = 1),
  \end{multline*}
which gives
  \begin{align*}
  P(Y_{j:n} - X_{j:n} > y, X_{j:n} \in (x,x+\ep))
    &\ge j\binom{n}{j}P(Y_{j:n} > x + \ep + y, A)\\
  &= j\binom{n}{j}P(\#\{i: Y_i < x + \ep + y\} < j, A)\\
  &\ge j\binom{n}{j}P(\#\{i \ne j: Y_i < x +\ep + y\} < j - 1, A).
  \end{align*}
Note that
  \begin{multline*}
  P(\#\{i \ne j: Y_i < x +\ep + y\} < j - 1, A) \\
    = \sum_{\ell=0}^{j-2}\sum_{m=\ell}^{j-2}
    P(\#\{i < j: Y_i < x + \ep + y\} = \ell,
    \#\{i > j: Y_i < x + \ep + y\} = m - \ell, A).
  \end{multline*}
If we now define $\wh q_i(x,y)=q_i(x,y+\ep)$, $\ol q_i(x,y)=q_i(x+\ep, y)$, $\wh p_i=1-\wh q_i$, and $\ol p_i=1-\ol q_i$, then this gives
  \begin{multline*}
  P(Y_{j:n} - X_{j:n} > y, X_{j:n} \in (x,x+\ep))\\
    \ge j\binom{n}{j}\sum_{\ell=0}^{j-2}\sum_{m=\ell}^{j-2}
    \binom{j-1}{\ell}\binom{n-j}{m-\ell}
    \wh p_1^\ell\wh q_1^{j-1-\ell}\ol p_2^{n-j-m+\ell}\ol q_2^{m-\ell}\\
    \cdot\Phi(x)^{j-1}\ol\Phi(x+\ep)^{n-j}P(X\in(x,x+\ep))
  \end{multline*}
We compare this with
  \begin{align*}
  \ph_{j:n}^>(x,y)
    &= P\bigg(\sum_{i=1}^{j-1} 1_{\{U_i\le q_1\}}
    > \sum_{i=j+1}^n 1_{\{U_i\le q_2\}}\bigg)\\
  &= \sum_{\ell=0}^{j-2}\sum_{m=\ell}^{j-2}
    P\bigg(\sum_{i=1}^{j-1} 1_{\{U_i\le q_1\}} = j - 1 - \ell,
    \sum_{i=j+1}^n 1_{\{U_i\le q_2\}} = m - \ell\bigg)\\
  &= \sum_{\ell=0}^{j-2}\sum_{m=\ell}^{j-2}
    \binom{j-1}{\ell}\binom{n-j}{m-\ell}
    p_1^\ell q_1^{j-1-\ell}p_2^{n-j-m+\ell}q_2^{m-\ell},
  \end{align*}
and the remainder of the proof is the same.

Using \eqref{RWrep.3}, we now prove the upper bound in \eqref{RWrep.2}. Let
  \[
  \xi = P(Y_{j:n} - X_{j:n} > y \mid X_{j:n}) - \ph^\ge_{j:n}(X_{j:n},y),
  \]
so that $\{\xi>0\}=\{X_{j:n}\in B\}$ for some Borel subset $B\subset\RR$. Fix $\ep>0$. There exists $B_0\subset\RR$ such that $B_0$ is a finite, disjoint union of open intervals, and
  \[
  |P(X_{j:n} \in B) - P(X_{j:n} \in B_0)| < \ep.
  \]
(See Proposition 1.20 in \cite{Fo99}, for example.) Hence, by \eqref{RWrep.3},
  \begin{align*}
  E[P(Y_{j:n} - X_{j:n} > y \mid X_{j:n})1_{\{X_{j:n} \in B\}}]
    &= P(Y_{j:n} - X_{j:n} > y, X_{j:n} \in B)\\
  &\le P(Y_{j:n} - X_{j:n} > y, X_{j:n} \in B_0) + \ep\\
  &\le E[\ph^\ge_{j:n}(X_{j:n},y)1_{\{X_{j:n} \in B_0\}}] + \ep\\
  &\le E[\ph^\ge_{j:n}(X_{j:n},y)1_{\{X_{j:n} \in B\}}] + 2\ep.
  \end{align*}
Therefore, $E[\xi1_{\{\xi>0\}}]\le 2\ep$. Letting $\ep\to0$ shows that $\xi \le0$ a.s., completing the proof. The lower bound in \eqref{RWrep.2} is proved similarly. \qed


\bigskip

With Theorem \ref{T:RWrep}, we have accomplished the first step of connecting our quantiles to a random walk. The second step, given by the next theorem, is to derive an estimate for this walk.

\begin{thm}\label{T:RWest}
Fix $\al\in(0,1)$. Let
  \[
  \begin{split}
  \si &= \si(r_1,r_2) = \al r_1 + (1 - \al)r_2,\\
  \mu &= \mu(r_1,r_2) = \al r_1 - (1 - \al)r_2,
  \end{split}
  \]
and suppose $j(n)/n=\al+o(n^{-1/2})$. Then for each $\tau>1$, there exist constants $C>0$ and $n_0\in\NN$ such that, for all $n\ge n_0$,
  \[
  \psi^\le_{j(n):n}(r_1,r_2)
    \le C\frac{\si^\tau}{n^\tau\mu^{2\tau}},
  \]
whenever $\mu>0$. Note that $C$ does not depend on $r_1$ or $r_2$.
\end{thm}

\bigskip

\pf
First note that if $n\al r_1\le 2$, then
  \[
  \frac{\si^\tau}{n^\tau\mu^{2\tau}}
    \ge \frac{(\al r_1)^\tau}{n^\tau(\al r_1)^{2\tau}}
    \ge 2^{-\tau} \ge 2^{-\tau}\psi^\le_{j(n):n}(r_1,r_2).
  \]
Hence, we may assume that $n\al r_1>2$.

Now, since $j(n)/n=\al+n^{-1/2}a_n$, where $a_n\to0$, we may write
  \[
  \psi^\le_{j(n):n}(r_1,r_2) = P(n^{-1}(\xi_L - \xi_U) + b_n \le -\mu)
    \le P(|n^{-1}(\xi_L - \xi_U) + b_n| \ge \mu),
  \]
where
  \begin{align*}
  \xi_L &= \sum_{i=1}^{j(n)-1} (1_{\{U_i\le r_1\}} - r_1),\\
  \xi_U &= \sum_{i=j(n)+1}^n (1_{\{U_i\le r_2\}} - r_2),\\
  b_n &= n^{-1/2}a_n(r_1 + r_2) - n^{-1}r_1.
  \end{align*}
By Chebyshev's inequality,
  \[
  \begin{split}
  \psi^\le_{j(n):n}(r_1,r_2)
    &\le \mu^{-2\tau}E|n^{-1}(\xi_L - \xi_U) + b_n|^{2\tau}\\
  &\le C(n\mu)^{-2\tau}(E|\xi_L|^{2\tau} + E|\xi_U|^{2\tau})
    + C\mu^{-2\tau}|b_n|^{2\tau}.
  \end{split}
  \]
Note that
  \[
  |b_n| \le Cn^{-1/2}(r_1 + r_2)
    \le \frac{Cn^{-1/2}\si}{\al\wedge(1 - \al)},
  \]
which implies $|b_n|^{2\tau}\le Cn^{-\tau}\si^{2\tau}$. It will therefore suffice to show that
  \[
  E|\xi_L|^{2\tau} + E|\xi_U|^{2\tau} \le C(n\si)^\tau.
  \]
By Lemma \ref{L:RWest}, there exists $n_0$ and $C$ such that $n\ge n_0$ implies
  \[
  E|\xi_L|^{2\tau} \le C((j(n)r_1)^\tau \vee (j(n)r_1)).
  \]
By making $n_0$ larger if necessary, we may assume $n\al/2\le j(n)\le 3n\al/2$. Since $n\al r_1>2$, this gives $j(n)r_1>1$, so that
  \[
  E|\xi_L|^{2\tau} \le C(j(n)r_1)^\tau
    \le C(n\al r_1)^\tau \le C(n\si)^\tau.
  \]
Similarly, for $n$ sufficiently large,
  \[
  \begin{split}
  E|\xi_U|^{2\tau}
    &\le C(((n-j(n))r_2)^\tau \vee ((n-j(n))r_2))\\
  &\le C((n(1 - \al)r_2)^\tau \vee 1)\\
  &\le C((n(1 - \al)r_2)^\tau \vee (n\al r_1)^\tau)
    \le C(n\si)^\tau,
  \end{split}
  \]
completing the proof. \qed

\section{Parameter estimates and gap regimes}\label{S:param_est}

Let us now return to the specific assumptions of our model, as outlined in the beginning of Section \ref{S:fdd}. Fix $0\le s<t$. We shall adopt the notation and definitions of Section \ref{S:RWrep}, in the special case that $X=B(s)-q(s)$ and $Y=B(t)-q(t)$. In this case,
  \begin{equation}\label{q_def_specific}
  \begin{split}
  q_1(x,y) &= P(B(t) > q(t) + x + y \mid B(s) < q(s) + x),\\
  q_2(x,y) &= P(B(t) < q(t) + x + y \mid B(s) > q(s) + x).
  \end{split}
  \end{equation}
Let us also define
  \begin{equation}\label{Psidef}
  \Psi(x,y) = P(B(t) > q(t) + x + y, B(s) < q(s) + x).
  \end{equation}
Our objective is to verify Condition (ii) of Theorem \ref{T:moment}. Ideally, we would like to show that
  \begin{equation}\label{tight_goal}
  P(|F_n(t) - F_n(s)| > \ep) \le C_p\ep^{-p}|t - s|^{p/4}.
  \end{equation}
In the end, we will actually show something slightly weaker, although our final estimate will be sufficient to verify the conditions of Theorem \ref{T:moment}. Our approach will begin by conditioning on $F_n(s)$, so that we may apply Theorem \ref{T:RWrep}. We will then use Theorem \ref{T:RWest} to obtain the specific bound we need. Implementing this strategy will require precise estimates on the function $q_1$ and $q_2$, in terms of $x$, $y$, $s$, and $t$. These estimates will come from Taylor expansions. We therefore begin with a Taylor expansion for $\Psi$.

\begin{prop}\label{P:Psi_Taylor}
Fix $T,K>0$. There exists a constant $C$ such that for all $0\le s<t\le T$ and all $x,y$ satisfying $|x|+|y|\le K$,
  \begin{equation}\label{Psi_Taylor}
  \Psi(x,y) = \Psi(0,0) - \frac12u(q(s),s)y
    + \frac1{2\sqrt{2\pi\de}}u(q(s),s)y^2 + R,
  \end{equation}
where
  \begin{equation}\label{Psi_Taylor_rem}
  |R| \le C\left({(|x| + |y|)(\de^{1/2} + |y| + \de^{-1/2}|y|^2)
    + \de^{-3/2}|y|^4}\right),
  \end{equation}
and $\de=t-s$.
\end{prop}

\pf Fix $T,K>0$, $s,t\in[0,T]$, and $x,y\in\RR$ such that $s<t$ and $|x|+|y| \le K$. Let $\de=t-s$. In what follows, $C$ shall denote a constant that depends only on $T$ and $K$, and may change value from line to line.

By Taylor's theorem, we may write
  \begin{multline}\label{Psi_Taylor1}
  \Psi(x,y) = \Psi(0,0) + \pa_x\Psi(0,0)x + \pa_y\Psi(0,0)y\\
    + \frac12\pa_x^2\Psi(0,0)x^2 + \pa_x\pa_y\Psi(0,0)xy
    + \frac12\pa_y^2\Psi(0,0)y^2\\
    + \frac16\pa_x^3\Psi(\th x,\th y)x^3
    + \frac12\pa_x^2\pa_y\Psi(\th x,\th y)x^2 y
    + \frac12\pa_x\pa_y^2\Psi(\th x,\th y)x y^2
    + \frac16\pa_y^3\Psi(\th x,\th y)y^3,
  \end{multline}
where $\th\in(0,1)$. Let $\Phi(x) =(2\pi)^{-1/2}\int_{-\infty}^x e^{-y^2/2}\,dy$ and $\ol\Phi=1-\Phi$. We first establish that for all integers $i\ge0$ and $j\ge1$,
  \begin{align}
  \pa_x^i\Psi &= \int_{-\infty}^{q(s)+x}\ol\Phi\left({
    \frac{x+y+q(t)-z}{\de^{1/2}}}\right)\pa_x^i u(z,s)\,dz,\label{dform1}\\
  \pa_x^i\pa_y^j\Psi &= -\de^{-(j-1)/2}
    \ol\Phi^{(j-1)}\left({\frac{y+q(t)-q(s)}{\de^{1/2}}}\right)
    \pa_x^i u(q(s)+x,s) + \pa_x^{i+1}\pa_y^{j-1}\Psi.\label{dform2}
  \end{align}
If $i=0$, then \eqref{dform1} follows directly from the definition of $\Psi$, \eqref{Psidef}. Differentiating \eqref{dform1} gives
  \begin{multline*}
  \pa_x^{i+1}\Psi = \ol\Phi\left({\frac{y+q(t)-q(s)}{\de^{1/2}}}\right)
    \pa_x^i u(q(s)+x,s)\\
    + \de^{-1/2}\int_{-\infty}^{q(s)+x}\ol\Phi'\left({
    \frac{x+y+q(t)-z}{\de^{1/2}}}\right)\pa_x^i u(z,s)\,dz.
  \end{multline*}
Applying integration by parts shows that \eqref{dform1} holds for $i+1$. By induction, this proves \eqref{dform1}. For \eqref{dform2}, let $j=1$ and let $i\ge0$ be arbitrary. By \eqref{dform1}, we have
  \begin{align*}
  \pa_x^i\pa_y\Psi &= \pa_y\left[{\int_{-\infty}^{q(s)+x}\ol\Phi\left({
    \frac{x+y+q(t)-z}{\de^{1/2}}}\right)\pa_x^i u(z,s)\,dz}\right]\\
  &= \pa_x\left[{\int_{-\infty}^{q(s)+x}\ol\Phi\left({
    \frac{x+y+q(t)-z}{\de^{1/2}}}\right)\pa_x^i u(z,s)\,dz}\right]\\
  &\quad - \ol\Phi\left({\frac{y+q(t)-q(s)}{\de^{1/2}}}\right)
    \pa_x^i u(q(s)+x,s)\\
  &= - \ol\Phi\left({\frac{y+q(t)-q(s)}{\de^{1/2}}}\right)
    \pa_x^i u(q(s)+x,s) + \pa_x^{i+1}\Psi,
  \end{align*}
which is \eqref{dform2}. Now assume \eqref{dform2} is valid for some $j\ge1$ and all $i\ge0$. Applying $\pa_y$ to both sides of \eqref{dform2} immediately shows that \eqref{dform2} is valid for $j+1$ and any $i\ge0$. By induction, this proves \eqref{dform2}.

Now, by \eqref{dform1}, we may write,
for any $i\ge 1$,
  \begin{align*}
  \pa_x^i\Psi(x,y) &= \int_{-\infty}^{q(s)+x}\pa_x^i u(z,s)\,dz
    - \int_{-\infty}^{q(s)+x}\int_{-\infty}^{q(t)+x+y}
    p(t-s,z,w)\pa_x^i u(z,s)\,dw\,dz\\
  &= \pa_x^{i-1}u(q(s)+x,s) - \ze(s,t,x,x+y),
  \end{align*}
where $\ze$ is given by \eqref{zeta_def}.
We shall adopt the convention that $\pa_x^{-1}u(x,t) := P(B(t)\le x)$, so that the above equality is, in fact, valid for all $i\ge 0$. Hence,
  \[
  \pa_x^i\Psi(x,y) = \pa_x^{i-1}u(q(s)+x,s)
    + \int_s^t \pa_s\ze(r,t,x,x+y)\,dr - \ze(t,t,x,x+y).
  \]
Since
  \[
  \ze(s,t,w,z) = \int_{-\infty}^{q(s)+w} P^x(B(t-s) \le q(t) + z)
    \pa_x^i u(x,s)\,dx,
  \]
it follows that
  \[
  \ze(t,t,w,z) = \int_{-\infty}^{q(t)+(w\wedge z)}\pa_x^i u(x,t)\,dx
    = \pa_x^{i-1} u(q(t)+(w\wedge z),t).
  \]
Therefore,
  \[
  \pa_x^i\Psi(x,y) = \pa_x^{i-1}u(q(s) + x,s)
    - \pa_x^{i-1}u(q(t) + x + (y\wedge0),s)
    + \int_s^t \pa_s\ze(r,t,x,x+y)\,dr.
  \]
By the mean value theorem, since $|q(t)-q(s)|\le C|t-s|=C\de$, we have
  \[
  |\pa_x^{i-1}u(q(s) + x,s) - \pa_x^{i-1}u(q(t) + x + (y\wedge0),s)|
    \le C(\de + |y|).
  \]
By \eqref{newIV.3}, $|\pa_s\ze(r,t,x,x+y)|\le C|t-r|^{-1/2}$, and so we obtain
  \begin{equation}\label{pa_x_est}
  |\pa_x^i\Psi(x,y)| \le C(\de^{1/2}+|y|),
  \end{equation}
for any $i\ge0$.
(Here we have used the fact that since $\de\le T$, there exists $C$ such that $\de\le C\de^{1/2}$ for all $\de\le T$.)

We next consider the derivatives with respect to $y$. By \eqref{dform2},
  \[
  \pa_y\Psi(0,0) = -\ol\Phi\left({\frac{q(t)-q(s)}{\de^{1/2}}}\right)
    u(q(s),s) + \pa_x\Psi(0,0).
  \]
By \eqref{pa_x_est}, $|\pa_x\Psi(0,0)|\le C\de^{1/2}$. Also, $|q(t)-q(s)| \le C\de$ and $\ol\Phi(x)=1/2+O(|x|)$. Hence,
  \[
  \left|{\pa_y\Psi(0,0) + \frac12 u(q(s),s)}\right| \le C\de^{1/2}.
  \]
Similarly,
  \[
  |\pa_x\pa_y\Psi(0,0)| = \left|{-\ol\Phi\left({
    \frac{q(t)-q(s)}{\de^{1/2}}}\right)\pa_x u(q(s),s)
    + \pa_x^2\Psi(0,0)}\right| \le C,
  \]
and
  \[
  \pa_y^2\Psi(0,0) = -\de^{-1/2}\ol\Phi'\left({
    \frac{q(t)-q(s)}{\de^{1/2}}}\right)u(q(s),s)
    + \pa_x\pa_y\Psi(0,0).
  \]
Since $\ol\Phi'(x)=-(2\pi)^{-1/2}e^{-x^2/2}=-(2\pi)^{-1/2}(1+O(|x|^2))$, this implies
  \[
  \left|{\pa_y^2\Psi(0,0) - \frac1{\sqrt{2\pi\de}}u(q(s),s)}\right|
    \le C.
  \]
For the third-order partial derivatives, we have
  \begin{align*}
  |\pa_x^2\pa_y\Psi(x,y)| &= \left|{-\ol\Phi\left({
    \frac{y+q(t)-q(s)}{\de^{1/2}}}\right)\pa_x^2 u(q(s)+x,s)
    + \pa_x^3\Psi(x,y)}\right| \le C,\\
  |\pa_x\pa_y^2\Psi(x,y)| &= \left|{-\de^{-1/2}\ol\Phi'\left({
    \frac{y+q(t)-q(s)}{\de^{1/2}}}\right)\pa_x u(q(s)+x,s)
    + \pa_x^2\pa_y\Psi(x,y)}\right| \le C\de^{-1/2},\\
  |\pa_y^3\Psi(x,y)| &= \left|{-\de^{-1}\ol\Phi''\left({
    \frac{y+q(t)-q(s)}{\de^{1/2}}}\right)u(q(s)+x,s)
    + \pa_x\pa_y^2\Psi(x,y)}\right|\\
    &\le C(\de^{-3/2}|y| + \de^{-1/2}),
  \end{align*}
where the last inequality follows from $\ol\Phi''(x)\le C|x|$.

Substituting all of this into \eqref{Psi_Taylor1} gives \eqref{Psi_Taylor}, where $R$ satisfies
  \begin{multline*}
  |R| \le C(\de^{1/2}|x| + \de^{1/2}|y| + \de^{1/2}|x|^2 + |xy|
    + |y|^2\\
  + \de^{1/2}|x|^3 + |x^3 y| + |x^2 y| + \de^{-1/2}|xy^2|
    + \de^{-3/2}|y|^4 + \de^{-1/2}|y|^3).
  \end{multline*}
Since $x$ and $y$ are restricted to a compact set, we may simplify this to
  \[
  |R| \le C(\de^{1/2}|x| + \de^{1/2}|y| + |xy|
    + |y|^2 + \de^{-1/2}|xy^2|
    + \de^{-3/2}|y|^4 + \de^{-1/2}|y|^3),
  \]
which is precisely \eqref{Psi_Taylor_rem}. \qed

\begin{corollary}\label{C:Psi_Taylor}
Recall $q_1,q_2$ given by \eqref{q_def_specific}. Fix $T,K>0$. There exist constants $\de_0$ and $C$ such that for all $0\le s<t\le T$ and all $x,y$ satisfying $|x|+|y|\le K$ and $|x|\le\de_0$,
  \begin{align}
  \al q_1(x,y) &= \Psi(0,0) - \frac12u(q(s),s)y
    + \frac1{2\sqrt{2\pi\de}}u(q(s),s)y^2 + R_1,\label{q1_Taylor}\\
  (1 - \al)q_2(x,y) &= \Psi(0,0) + \frac12u(q(s),s)y
    + \frac1{2\sqrt{2\pi\de}}u(q(s),s)y^2 + R_2,\label{q2_Taylor}
  \end{align}
where $R_1,R_2$ both satisfy \eqref{Psi_Taylor_rem}, and $\de=t-s$.
\end{corollary}

\pf By \eqref{Psi_Taylor} and \eqref{pa_x_est},
  \begin{equation}\label{Psi_bound}
  |\Psi(x,y)| \le C(\de^{1/2} + |y| + \de^{-1/2}|y|^2) + |R|,
  \end{equation}
where $R$ satisfies \eqref{Psi_Taylor_rem}. By \eqref{q_def_specific} and \eqref{Psidef},
  \[
  q_1(x,y) = \frac{\Psi(x,y)}{P(B(s) < q(s) + x)}.
  \]
Note that $P(B(s)<q(s)+x)=\al+r(x,s)$, where $|r(x,s)|\le C|x|$. Hence,
  \[
  |\al q_1(x,y) - \Psi(x,y)| = \left|{
    \frac{r(x,s)}{\al + r(x,s)}}\right|\Psi(x,y)
    \le \frac{C|x|}{\al - C|x|}\Psi(x,y).
  \]
If $x$ is sufficiently small, then using \eqref{Psi_bound}, we have
  \[
  |\al q_1(x,y) - \Psi(x,y)|
    \le C(\de^{1/2}|x| + |xy| + \de^{-1/2}|xy^2|) + |R|.
  \]
By \eqref{Psi_Taylor}, this gives
  \begin{multline*}
  \bigg|\al q_1(x,y) - \Psi(0,0) + \frac12u(q(s),s)y
    - \frac1{2\sqrt{2\pi\de}}u(q(s),s)y^2\bigg|\\
    \le C(\de^{1/2}|x| + |xy| + \de^{-1/2}|xy^2|) + 2|R|.
  \end{multline*}
Since this is bounded above by the right-hand side of \eqref{Psi_Taylor_rem}, this completes the proof of \eqref{q1_Taylor}.

For \eqref{q2_Taylor}, let $\wt B(t)=-B(t)$, and let $\wt q(t)$ be the $(1-\al)$-quantile of the law of $\wt B(t)$, so that $\wt q(t)=-q(t)$. Let $\wt u(x,t)$ be the density of $\wt B(t)$, so that $\wt u(x,t)=u(-x,t)$. Define
  \[
  \wt\Psi(x,y) = P(\wt B(t) > \wt q(t) + x + y,
    \wt B(s) < \wt q(s) + x),
  \]
and
  \begin{align*}
  \wt q_1(x,y) &= P(\wt B(t) > \wt q(t) + x + y \mid
    \wt B(s) < \wt q(s) + x)\\
  &= P(B(t) < q(t) - x - y \mid B(s) > q(s) - x)\\
  &= q_2(-x,-y).
  \end{align*}
Hence, by \eqref{q1_Taylor},
  \begin{align*}
  (1-\al)q_2(x,y) &= (1-\al)\wt q_1(-x,-y)\\
  &= \wt\Psi(0,0) + \frac12\wt u(\wt q(s),s)y
    + \frac1{\sqrt{2\pi\de}}\wt u(\wt q(s),s)y^2 + R_2\\
  &= \wt\Psi(0,0) + \frac12u(q(s),s)y
    + \frac1{\sqrt{2\pi\de}}u(q(s),s)y^2 + R_2,
  \end{align*}
where $R_2$ satisfies \eqref{Psi_Taylor_rem}. To complete the proof, we observe that
  \begin{align*}
  \wt\Psi(0,0) &= P(B(t) < q(t), B(s) > q(s))\\
  &= P(B(t) < q(t)) - P(B(t) < q(t), B(s) < q(s))\\
  &= P(B(s) < q(s)) - P(B(t) < q(t), B(s) < q(s))\\
  &= P(B(t) > q(t), B(s) < q(s)) = \Psi(0,0),
  \end{align*}
which gives \eqref{q2_Taylor}. \qed

\bigskip

We are now ready to establish \eqref{tight_goal} and complete the proof of our main result, Theorem \ref{T:main}. Recall that $F_n=n^{1/2}(Q_n-q)$. Hence, the event $\{|F_n(t)-F_n(s)|>\ep\}$ is precisely the event that the process $Q_n-q$ traverses a gap of size $n^{-1/2}\ep$ in the time interval $[s,t]$. In that same time interval, each Brownian particle will traverse a gap of size roughly $|t-s|^{1/2}$. We therefore divide our proof into three separate cases, or parameter regimes: $n^{-1/2}\ep\gg|t-s|^{1/2}$ (the large gap regime), $n^{-1/2} \ep\ll |t-s|^{1/2}$ (the small gap regime), and $n^{-1/2}\ep\approx|t-s|^{1/2}$ (the medium gap regime).

The proof in the large gap regime (Lemma \ref{L:large_gap}) uses only crude estimates and requires none of the machinery we have so far developed. The proof is entirely standard and is given in the appendix. The proof in the small gap regime (Lemma \ref{L:small_gap}) requires the most care. We follow the strategy outlined earlier, utilizing all of the tools developed previously. The proof in the medium gap regime (Lemma \ref{L:medium_gap}) is a minor modification of the proof of Lemma \ref{L:small_gap}. In making this modification, however, our estimate loses some precision, and the result in the medium gap regime is, in fact, slightly weaker than \eqref{tight_goal}.


\begin{lemma}\label{L:large_gap}
Fix $T>0$, $\De\in(0,1/2)$, and $p>2$. There exists a positive constant $C$ such that
  \[
  P(|F_n(t) - F_n(s)| > \ep) \le C(\ep^{-1}|t - s|^{1/4})^p,
  \]
whenever $s,t\in[0,T]$, $0<\ep<1$, and $n^{-1/2}\ep\ge |t-s|^{1/2-\De}$.
\end{lemma}

\begin{lemma}\label{L:small_gap}
Fix $T>0$, $\De\in(0,1/2)$, and $p>2$. There exists a positive constant $C$ and an integer $n_0$ such that
  \[
  P(|F_n(t) - F_n(s)| > \ep) \le C(\ep^{-1}|t - s|^{1/4})^p,
  \]
whenever $n\ge n_0$, $s,t\in[0,T]$, $0<\ep<1$, and $n^{-1/2}\ep\le |t-s|^{1/2+\De}$.
\end{lemma}

\pf Fix $T>0$, $\De\in(0,1/2)$, and $p>2$. Without loss of generality, assume $s<t$ and define $\de=t-s$. Note that it is sufficient to prove that there exists a constant $\de_0>0$ such that the lemma holds whenever $\de\le\de_0$.

Let $n_0$ be as in Theorem \ref{T:RWest}, with $\tau=p/2$. Define $\ol B_j=B_j-q$ and $\ol Q_n=Q_n-q$. Note that $\ol Q_n=\ol B_{j(n):n}$ and $F_n=n^{1/2}\ol Q_n$. Let $K=n^{-1/2}\ep\de^{-1/4}$ and $y=n^{-1/2} \ep$. Then by Theorem \ref{T:RWrep},
  \begin{align}
  P(F_n(t) - F_n(s) < -\ep) &= P(\ol Q_n(t) - \ol Q_n(s) < -y)\notag\\
  &= E[P(\ol Q_n(t) - \ol Q_n(s) < -y \mid \ol Q_n(s))]\notag\\
  &\le E[\ph_{j(n):n}^\le(\ol Q_n(s),-y)]\notag\\
  &\le \sup_{|x|\le K}\ph_{j(n):n}^\le(x,-y) + P(|\ol Q_n(s)|>K).
    \label{small_gap1}
  \end{align}
By making $n_0$ larger, if necessary, we may apply Proposition \ref{P:quant_tails}
and Remark \ref{R:quant_tails}
to obtain
  \begin{equation}\label{small_gap2}
  P(|\ol Q_n(s)| > K) = P(|F_n(s)| > \ep\de^{-1/4})
    \le C(\ep\de^{-1/4})^{-p}
    = C(\ep^{-1}|t - s|^{1/4})^p.
  \end{equation}
For the other term, fix $x\in[-K,K]$. Following the notation of Theorem \ref{T:RWest}, define
  \begin{align*}
  \si(x,y) &= \al q_1(x,y) + (1 - \al)q_2(x,y),\\
  \mu(x,y) &= \al q_1(x,y) - (1 - \al)q_2(x,y).
  \end{align*}
By Corollary \ref{C:Psi_Taylor},
  \[
  \mu(x,-y) = u(q(s),s)y + R_\mu,
  \]
where $R_\mu$ satisfies \eqref{Psi_Taylor_rem}. Note that $|x|\le K=n^{-1/2} \ep\de^{-1/4}$ and $y=n^{-1/2}\ep\le\de^{1/2+\De}$. In particular, $y\le K$ and $y\le\de^{1/2}$. Therefore, \eqref{Psi_Taylor_rem} simplifies to
  \[
  |R_\mu| \le C(K\de^{1/2} + \de^{-3/2}y^4)
    \le C(K\de^{1/2} + \de^{3\De}y)
    = C(\de^{1/4} + \de^{3\De})y.
  \]
Hence, since $s\mapsto u(q(s),s)$ is bounded below on $[0,T]$, there exists $\de_0$ such that $\de\le\de_0$ implies $\mu(x,-y) > Cy$. Similarly,
  \[
  \si(x,-y) = 2\Psi(0,0)
    + \frac1{\sqrt{2\pi\de}}u(q(s),s)y^2 + R_\si,
  \]
where $|R_\si|\le C(\de^{1/4} + \de^{3\De})y\le Cy\le C\de^{1/2}$. By \eqref{pa_x_est}, this implies $|\si(x,-y)|\le C\de^{1/2}$.

We now apply Theorem \ref{T:RWest} with $\tau=p/2$ to obtain
  \begin{equation}\label{small_gap3}
  \ph_{j(n):n}^\le(x,-y)
    \le C\left({\frac{\si(x,-y)}{n\mu(x,-y)^2}
      }\right)^{p/2}
    \le C\left({\frac{\de^{1/2}}{ny^2}
      }\right)^{p/2}
    = C(\ep^{-1}|t - s|^{1/4})^p.
  \end{equation}
Substituting \eqref{small_gap3} and \eqref{small_gap2} into \eqref{small_gap1} gives
  \[
  P(F_n(t) - F_n(s) < -\ep) \le C(\ep^{-1}|t - s|^{1/4})^p.
  \]
To complete the proof, we follow as in \eqref{small_gap1} to obtain
  \[
  P(F_n(t) - F_n(s) > \ep)
    \le \sup_{|x|\le K}\ph_{j(n):n}^\ge(x,y)
    + P(|\ol Q_n(s)| > K),
  \]
and apply \eqref{small_gap2} to the second term. For the first term, recall from Section \ref{S:RWrep} that
  \[
  \ph_{j(n):n}^\ge(x,y)
    = \psi_{(n-j(n)+1):n}^\le(q_2(x,y),q_1(x,y)).
  \]
Since $(n-j(n)+1)/n = (1-\al) + o(n^{-1/2})$, we may apply Theorem \ref{T:RWest} to obtain
  \[
  \ph_{j(n):n}^\ge(x,y)
    \le C\left({\frac{\wt\si(x,y)}{n\wt\mu(x,y)^2}}\right)^{p/2},
  \]
where
  \begin{align*}
  \wt\si(x,y) &= (1 - \al)q_2(x,y) + \al q_1(x,y) = \si(x,y),\\
  \wt\mu(x,y) &= (1 - \al)q_2(x,y) - \al q_1(x,y) = -\mu(x,y).
  \end{align*}
The estimates preceding \eqref{small_gap3} can now be used, as before, to give
  \[
  \ph_{j(n):n}^\ge(x,y) \le C(\ep^{-1}|t - s|^{1/4})^p,
  \]
which completes the proof. \qed


%
\begin{lemma}\label{L:medium_gap}
Fix $T>0$, $\De\in(0,1/8)$, and $p>2$. There exists a positive constant $C$ and an integer $n_0$ such that
  \[
  P(|F_n(t) - F_n(s)| > \ep)
    \le C(\ep^{-1}|t - s|^{1/4-2\De})^p,
  \]
whenever $n\ge n_0$, $s,t\in[0,T]$, $0<\ep<1$, and $|t-s|^{1/2+\De} \le n^{-1/2}\ep\le |t-s|^{1/2-\De}$.
\end{lemma}

\pf Fix $T>0$, $\De\in(0,1/8)$, and $p>2$. Without loss of generality, assume $s<t$ and define $\de=t-s$. Note that it is sufficient to prove that there exists a constant $\de_0>0$ such that the lemma holds whenever $\de\le\de_0$.

Let $n_0$ be as in Theorem \ref{T:RWest}, with $\tau=p/2$. Let $K=n^{-1/2} \ep\de^{-1/4}$ and $y=n^{-1/2} \ep$. Also define $\wt\ep=n^{1/2}\de^{1/2 +\De}\le\ep$ and $\wt y=n^{-1/2}\wt\ep=\de^{1/2+\De}$. As in \eqref{small_gap1} and \eqref{small_gap2},
  \begin{align*}
  P(F_n(t) - F_n(s) < -\ep) &\le P(F_n(t) - F_n(s) < -\wt\ep)\\
  &\le \sup_{|x|\le K}\ph_{j(n):n}^\le(x,-\wt y)
    + C(\ep^{-1}|t-s|^{1/4})^p.
  \end{align*}
By the estimates following \eqref{small_gap2}, $\mu(x,-\wt y)=u(q(s),s) \wt y + R_\mu$, where
  \[
  |R_\mu| \le C(K\de^{1/2} + \de^{-3/2}\wt y^4)
    \le C(\de^{3/4-\De} + \de^{1/2+4\De})
    = C(\de^{1/4-2\De} + \de^{3\De})\wt y.
  \]
Since $\De<1/8$ and $s\mapsto u(q(s),s)$ is bounded below on $[0,T]$, there exists $\de_0$ such that $\de\le\de_0$ implies $\mu(x,-\wt y) > C\wt y$. Similarly, $|\si(x,-\wt y)|\le C\de^{1/2}$. As in \eqref{small_gap3}, this gives $\ph_{j(n):n}^\le(x,-\wt y) \le C(n^{-1}\wt y^{-2}\de^{1/2})^{p/2}$. Note that $\wt y=\de^{1/2+\De} \ge n^{-1/2} \ep\de^{2\De}$. Hence,
  \[
  \ph_{j(n):n}^\le(x,-\wt y) \le C(\ep^{-2}\de^{1/2-4\De})^{p/2}
    = C(\ep^{-1}|t - s|^{1/4-2\De})^p.
  \]
As in the second half of the proof of Lemma \ref{L:small_gap}, the bound on $P(F_n(t)-F_n(s)>\ep)$ is proved similarly. \qed

\bigskip

\noindent{\bf Proof of Theorem \ref{T:main}.} By Corollary \ref{C:quant_CLT}, it will suffice to show that $\{F_n\}$ is relatively compact in $C[0,\infty)$. By Theorem \ref{T:moment}, 
we need only verify that
  \[
  \sup_{n\ge n_0} P(|F_n(t) - F_n(s)| \ge \ep)
    \le C\ep^{-\al}|t - s|^{1+\beta}.
  \]
Taking $\De=1/16$ and $p=9$ in Lemmas \ref{L:large_gap}, \ref{L:small_gap}, and \ref{L:medium_gap} shows that
  \[
  P(|F_n(t) - F_n(s)| \ge \ep) \le C\ep^{-9}|t - s|^{9/8},
  \]
for all $n\ge n_0$, which completes the proof. \qed

\section*{Acknowledgments}

The author is grateful for the efforts of three very thorough referees, whose comments and suggestions have greatly improved this manuscript.

\appendix

\section{Appendix}

In this section, we have collected together the more routine and straightforward proofs, so as not to distract from the essential techniques of the main paper. We begin with the proof of Lemma \ref{L:quant_ODE}.

\bigskip

\noindent{\bf Proof of Lemma \ref{L:quant_ODE}.} Standard results imply that for each $t_0>0$, there exists a neighborhood $I$ of $t_0$ on which \eqref{quant_ODE} has a unique solution $\wt q$ with $\wt q(t_0)=q(t_0)$. To show that $\wt q=q$ on $I$, it suffices to show that $g(t)=P(B(t)\le\wt q(t))$ satisfies $g=\al$ on $I$. Since $g(t_0)=\al$, it suffices to show $g'=0$ on $I$. For this, we simply differentiate
  \[
  g(t) = \int_{-\infty}^{\wt q(t)} u(x,t)\,dx,
  \]
which gives
  \[
  g'(t) = u(\wt q(t),t)\wt q'(t)
    + \int_{-\infty}^{\wt q(t)}\pa_t u(x,t)\,dx
    = -\frac12\pa_x u(\wt q(t),t)
    + \frac12\int_{-\infty}^{\wt q(t)}\pa_x^2 u(x,t)\,dx
    = 0.
  \]
Since $t_0>0$ was arbitrary, this shows that $q\in C^\infty(0,\infty)$ and satisfies \eqref{quant_ODE} for all $t>0$.

To show that $q(t)\to q(0)$ as $t\to0$, we first show that $q(0)\le\liminf_{t\to0} q(t)$. The proof is by contradiction. Suppose not. Then there exists $\ep>0$ and $t_n\downarrow0$ such that $q(t_n)<q(0)-\ep$ for all $n$. Hence,
  \begin{multline*}
  \al = P(B(t_n) \le q(t_n)) \le P(B(t_n) \le q(0) - \ep)\\
    \to P(B(0) \le q(0) - \ep) \le P(B(0) \le q(0)) = \al.
  \end{multline*}
It follows that $P(B(0)\le q(0)-\ep)=\al$, but this contradicts the uniqueness of the $\al$-quantile of the measure $f(x)\,dx$. The proof that $\limsup_{t\to0} q(t)\le q(0)$ is similar. \qed

\bigskip

We next present a proof of the quantile central limit theorem. The proof uses a multi-dimensional version of the Lindeberg-Feller theorem, which is stated below as Theorem \ref{T:Lind-Feller}.

\begin{thm}\label{T:Lind-Feller}
For each fixed $n$, let $\{X_{m,n}\}_{m=1}^n$ be independent, $\RR^d$-valued random vectors with mean zero and covariance matrix $\si_{m,n}$. If
  \begin{enumerate}[(i)]
  \item $\sum_{m=1}^n \si_{m,n} \to \sigma$ as $n\to\infty$, and
  \item for each $\th\in\RR^d$ and each $\ep>0$, $\sum_{m=1}^n E[|\th\cdot X_{m,n}|^2 1_{\{|\th\cdot X_{m,n}|>\ep\}}] \to 0$ as $n\to\infty$,
  \end{enumerate}
then $S_n=X_{1,n}+\cdots+X_{n,n}\To N$, where $N$ is multi-normal with mean $0$ and covariance $\si$.
\end{thm}

\noindent{\textbf{Proof of Theorem \ref{T:quant_CLT}.} }
Given $x,y\in\RR^d$, we shall write $x\le y$ if $x(i)\le y(i)$ for all $i$.

Fix $x\in\RR^d$ and for each $n,m\in\NN$, $1\le m\le n$, define the random vector $Y_{m,n}\in\RR^d$ by
  \[
  Y_{m,n}(j) = n^{-1/2}\left(
    1_{\{X_m(j) \le n^{-1/2}x(j) + q(j)\}} - p_n(j)\right),
  \]
where $p_n(j)=\Phi_j(n^{-1/2}x(j)+q(j))$. Then for each fixed $n\in\NN$, $Y_{m,n}$, $1\le m\le n$, are independent, $EY_{m,n}=0$, and
  \begin{enumerate}[(i)]
  \item $\sum_{m=1}^n E[Y_{m,n}(i)Y_{m,n}(j)]\to G_{ij}(q(i),q(j)) - \al(i)\al(j)$ as $n\to\infty$,
  \item for each $\th\in\RR^d$ and $\ep>0$, $\sum_{m=1}^n E[|\th\cdot Y_{m,n}|^2 1_{\{|\th\cdot Y_{m,n}|>\ep\}}] \to 0$ as $n\to\infty$.
  \end{enumerate}
Part (i) follows since
  \begin{multline*}
  \sum_{m=1}^n E[Y_{m,n}(i)Y_{m,n}(j)]\\
  = n^{-1}\sum_{m=1}^n [P(
    X_m(i) \le n^{-1/2}x(i) + q(i),
    X_m(j) \le n^{-1/2}x(j) + q(j))
    - p_n(i)p_n(j)]\\
  = P(X(i) \le n^{-1/2}x(i) + q(i),
    X(j) \le n^{-1/2}x(j) + q(j))
    - p_n(i)p_n(j),
  \end{multline*}
and part (ii) follows since $|\th\cdot Y_{m,n}|\le n^{-1/2}
d
\max(|\th(1)|, \ldots,|\th(d)|)$, and therefore $P(|\th\cdot Y_{m,n}| > \ep) = 0$ for sufficiently large $n$.

Thus, by Theorem \ref{T:Lind-Feller}, $S_n=Y_{1,n}+\cdots+Y_{n,n}\To\wt N$, where $\wt N$ is multinormal with mean $0$ and covariance
  \begin{equation}\label{quant_CLT.1}
  E[\wt N(i)\wt N(j)] = G_{ij}(q(i),q(j)) - \al(i)\al(j).
  \end{equation}
Now, observe that the following are equivalent:
  \begin{enumerate}[(a)]
  \item $n^{1/2}(X_{\ka(n):n}-q)\le x$,
  \item $X_{\ka(n):n}(j)\le n^{-1/2}x(j) + q(j)$ for all $j$,
  \item $\sum_{m=1}^n 1_{\{X_m(j) \le n^{-1/2}x(j) + q(j)\}}
    \ge \ka(n,j)$ for all $j$,
  \item $n^{-1/2}\sum_{m=1}^n(1_{\{X_m(j) \le n^{-1/2}x(j) + q(j)\}}
    - p_n(j)) \ge n^{-1/2}(\ka(n,j) - n p_n(j))$ for all $j$.
  \end{enumerate}
Thus, if $a_n\in\RR^d$ is defined by $a_n(j)=n^{-1/2}(\ka(n,j) - n p_n(j))$, then $P(n^{1/2}(X_{\ka(n):n}-q)\le x)=P(S_n\ge a_n)$. Note that
  \[
  a_n(j) = n^{1/2}(\al(j) + o(n^{-1/2}) - p_n(j))
    = \frac{\Phi_j(q(j)) - \Phi_j(n^{-1/2}x(j) + q(j))}{n^{-1/2}}
    + o(1),
  \]
so that $a_n\to a\in\RR^d$, where $a(j)=-x(j)\Phi_j'(q(j))$. Therefore,
  \[
  P(n^{1/2}(X_{\ka(n):n}-q)\le x) \to P(\wt N \ge a)
    = P(\wt N \le -a) = P(N\le x),
  \]
where $N$ is the random vector defined by $N(j)=\wt N(j)/\Phi_j'(q(j))$. Comparing \eqref{quant_CLT.1} and \eqref{quant_CLT}, this completes the proof. \qed

\bigskip

At certain points in Section \ref{S:param_est}, we use the fact that $F_n$ has well-behaved tail probabilities. The precise formulation of this fact is given below.

\begin{prop}\label{P:quant_tails}
Fix $T>0$. For each $r>0$, there exist constants $C$ and $n_0$ such that
  \[
  P(|F_n(t)| > \la) \le C(\la^{-r} + P(|B(t) - q(t)|^2 > \la)),
  \]
for all $n\ge n_0$, $0\le t\le T$, and $\la>0$. In particular, $\sup_n P(|F_n(0)|>\la) \to 0$ as $\la\to\infty$.
\end{prop}

\begin{remark}\label{R:quant_tails}
Since $B$ is Gaussian, Proposition \ref{P:quant_tails} in fact shows that $P(|F_n(t)| > \la) \le C\la^{-r}$. Indeed,
  \[
  P(|B(t) - q(t)|^2 > \la) \le P(|B(t)| + M > \la^{1/2})
    = 2P(B(t) \le -\la^{1/2} + M),
  \]
where $M=\sup_{0\le t\le T}|q(t)|$. If $\la>4M^2$, then
  \[
  P(|B(t) - q(t)|^2 > \la) \le 2P(B(t) \le -\la^{1/2}/2)
    = 2\Phi(-t^{-1/2}\la^{1/2}/2)
    \le 2\Phi(-T^{-1/2}\la^{1/2}/2),
  \]
where $\Phi$ is the distribution function of the standard normal, which satisfies $\Phi(-x)\le C_rx^{-r}$ for all $r>0$. Hence, $P(|B(t) - q(t)|^2 > \la) \le 2C_{2r}T^{-r}2^{-2r}\la^r$ for $\la$ sufficiently large.
\end{remark}

\noindent{\textbf{Proof of Proposition \ref{P:quant_tails}.}}
Fix $T>0$ and $r>0$. Since $u(q(0),0)>0$, there exists $\ep\in(0,1)$ such that
  \[
  m = \inf\{u(x,t): |x - q(t)| \le \ep, 0 \le t \le T\} > 0.
  \]
Since $j(n)/n=\al+o(n^{-1/2})$, we may choose $n_0 \ge 2m^{-1}$ such that $|j(n)/n-\al|\le n^{-1/2}\ep m/2$, for all $n\ge n_0$. Let $n\ge n_0$, $t\in[0,T]$, and $\la>0$ be arbitrary. Note that we may assume without loss of generality that $\la > 1+2m^{-1/2}$ and $r>2$.

Let us first assume $n^{-1/2}\la\le\ep$. We begin with
  \begin{align}
  P(F_n(t) < -\la)
    &= P\bigg(\sum_{j=1}^n 1_{\{B_j(t) < q(t) - n^{-1/2}\la\}}
    \ge j(n)\bigg)\notag\\
  &= P\bigg(n^{-1}\sum_{j=1}^n(1_{\{B_j(t) < q(t) - n^{-1/2}\la\}} - p_n)
    \ge \frac{j(n)}n - p_n\bigg),\label{quant_tails.1}
  \end{align}
where $p_n=p_n(\la)=P(B(t) < q(t) - n^{-1/2}\la)$. By the mean value theorem, $\al - p_n\ge n^{-1/2}\la m$. Since $|j(n)/n-\al|\le n^{-1/2}m/2\le n^{-1/2}\la m/2$, we have
  \[
  P(F_n(t) < -\la)
    \le P\bigg(n^{-1}\sum_{j=1}^n
    (1_{\{B_j(t) < q(t) - n^{-1/2}\la\}} - p_n)
    \ge n^{-1/2}\la m/2\bigg).
  \]
By Chebyshev's inequality,
  \[
  P(F_n(t) < -\la)
    \le Cn^{-r/2}\la^{-r}E\bigg|\sum_{j=1}^n
    (1_{\{B_j(t) < q(t) - n^{-1/2}\la\}} - p_n)\bigg|^r.
  \]
By Burkholder's inequality (see, for example, Theorem 6.3.10 in \cite{St}),
  \[
  P(F_n(t) < -\la)
    \le Cn^{-r/2}\la^{-r}E\bigg|\sum_{j=1}^n
    |1_{\{B_j(t) < q(t) - n^{-1/2}\la\}} - p_n|^2\bigg|^{r/2}.
  \]
Finally, by Jensen's inequality,
  \[
  P(F_n(t) < -\la)
    \le C\la^{-r}E|1_{\{B(t) < q(t) - n^{-1/2}\la\}} - p_n|^r.
  \]
Note that
  \begin{align*}
  E|1_{\{B(t) < q(t) - n^{-1/2}\la\}} - p_n|^r
    &= p_n(1 - p_n)^r + (1 - p_n)p_n^r\\
    &= p_n(1 - p_n)((1 - p_n)^{r-1} + p_n^{r-1}) \le 2p_n.
  \end{align*}
Hence, $P(F_n(t) < -\la)\le C\la^{-r}p_n(\la)\le C\la^{-r}$.

Similarly, using the mean value theorem,
  \begin{align*}
  P(F_n(t) > \la)
    &= P\bigg(\sum_{j=1}^n 1_{\{B_j(t) > q(t) + n^{-1/2}\la\}}
    \ge n - j(n) + 1\bigg)\\
  &\le P\bigg(n^{-1}\sum_{j=1}^n
    (1_{\{B_j(t) > q(t) + n^{-1/2}\la\}} - \ol p_n)
    \ge 1 - \frac{j(n)}n - \ol p_n\bigg)\\
  &\le P\bigg(n^{-1}\sum_{j=1}^n
    (1_{\{B_j(t) > q(t) + n^{-1/2}\la\}} - \ol p_n)
    \ge n^{-1/2}\la m/2\bigg).
  \end{align*}
where $\ol p_n=\ol p_n(\la)=P(B(t) > q(t) + n^{-1/2}\la)$. By Chebyshev, Burkholder, and Jensen, $P(F_n(t) > \la)\le C\la^{-r}\ol p_n(\la)\le C\la^{-r}$. This completes the proof in the case $n^{-1/2}\la\le\ep$.

We now assume $n^{-1/2}\la>\ep$. In this case, $\al-p_n\ge P(q(t)-\ep<B(t) <q(t))\ge\ep m$. Hence, $j(n)/n-p_n\ge \ep m/2$. Thus, starting from \eqref{quant_tails.1}, Chebyshev, Burkholder, and Jensen imply
  \[
  P(F_n(t) < -\la)\le Cn^{-2r}E\bigg|
    \sum_{j=1}^n(1_{\{B_j(t) < q(t) - n^{-1/2}\la\}} - p_n)\bigg|^{2r}
    \le Cn^{-r}p_n(\la).
  \]
Similarly, $P(F_n(t)>\la)\le Cn^{-r}\ol p_n(\la)$. Hence,
  \[
  P(|F_n(t)| > \la) \le Cn^{-r}(p_n(\la) + \ol p_n(\la))
    = Cn^{-r}P(|B(t) - q(t)| > n^{-1/2}\la).
  \]
If $n>\la$, then $P(|F_n(t)|>\la)\le C\la^{-r}$, and we are done. If $n\le \la$, then
  \[
  P(|F_n(t)| > \la) \le CP(|B(t) - q(t)| > \la^{1/2})
    = CP(|B(t) - q(t)|^2 > \la),
  \]
and this completes the proof. \qed

\bigskip

In the proof of Theorem \ref{T:RWrep}, we need the following formula for the distribution of the $j$-th order statistic.

\begin{lemma}\label{L:RWrep}
If $\{X_n\}$ is a sequence of iid copies of $X$, then for all $x\in \RR$,
  \[
  P(X_{j:n} \le x) = \int_{-\infty}^x j\binom{n}{j}
    \Phi(y)^{j-1}\ol\Phi(y)^{n-j}\,d\Phi(y),
  \]
where $\Phi(x)=P(X\le x)$ and $\ol\Phi=1-\Phi$.
\end{lemma}

\pf Recall that $\Phi$ is continuous. Hence, if $g$ is absolutely continuous on $[0,1]$, then
  \[
  g(\Phi(x)) = g(0) + \int_{-\infty}^x g'(\Phi(y))\,d\Phi(y).
  \]
(See Exercise 3.36(b) in \cite{Fo99}, for example.) We therefore have
  \[
  \begin{split}
  P(X_{j:n} \le x) &= P\bigg(\sum_{i=1}^n 1_{\{X_i\le x\}} \ge j\bigg)
    = \sum_{k=j}^n P\bigg(\sum_{i=1}^n 1_{\{X_i\le x\}} = k\bigg)
    = \sum_{k=j}^n \binom{n}{k}\Phi(x)^k\ol\Phi(x)^{n-k}\\
  &= \sum_{k=j}^n \binom{n}{k}\int_{-\infty}^x
    (k\Phi(y)^{k-1}\ol\Phi(y)^{n-k}
    - (n - k)\Phi(y)^k\ol\Phi(y)^{n-k-1})\,d\Phi(y).
  \end{split}
  \]
We can rewrite this as
  \begin{multline*}
  P(X_{j:n} \le x) = \int_{-\infty}^x\bigg[
    \sum_{k=j}^n k\binom{n}{k}\Phi(y)^{k-1}\ol\Phi(y)^{n-k}\\
    - \sum_{k=j+1}^n(n - k + 1)\binom{n}{k-1}
    \Phi(y)^{k-1}\ol\Phi(y)^{n-k}\bigg]\,d\Phi(y).
  \end{multline*}
Observing that
  \[
  (n - k + 1)\binom{n}{k-1} = k\binom{n}{k},
  \]
completes the proof. \qed

\bigskip

The following lemma is needed in the proof of Theorem \ref{T:RWest}.

\begin{lemma}\label{L:RWest}
For each $r\ge1$, there exist constants $C>0$ and $n_0\in\NN$ such that
  \begin{equation}\label{RWest}
  E\bigg|\sum_{i=1}^n(1_{\{U_i \le p\}} - p)\bigg|^{2r}
    \le C((np)^r \vee (np)),
  \end{equation}
for all $n\ge n_0$ and all $p\in[0,1]$. Note that $C$ depends only on $r$, and not on $n$ or $p$.
\end{lemma}

\begin{remark}
If $r<1$, then Lemma \ref{L:RWest}, together with Jensen's inequality, imply that the expectation in \eqref{RWest} is bounded above by $C(np)^r$.
\end{remark}

\noindent{\textbf{Proof of Lemma \ref{L:RWest}.} }
First assume that $r\in\NN$. Let $\xi=1_{\{ U_i\le p\}}-p$ and $\ph(t)=E[e^{it\xi}]$. Then
  \[
  E\bigg|\sum_{i=1}^n(1_{\{U_i \le p\}} - p)\bigg|^{2r}
    = (-1)^r\frac{d^{2r}}{dt^{2r}}(\ph(t)^n)\bigg|_{t=0}.
  \]
By Fa\`a di Bruno's formula (see \cite{Jo}, for example), if $n > 2r$, then
  \[
  \frac{d^{2r}}{dt^{2r}}(\ph(t)^n)
    = \sum \frac{(2r)!}{b_1!b_2!\cdots b_{2r}!}
    \frac{n!}{(n-k)!}\ph(t)^{n-k}
    \left({\frac{\ph'(t)}{1!}}\right)^{b_1}
    \left({\frac{\ph''(t)}{2!}}\right)^{b_2}\cdots
    \left({\frac{\ph^{(2r)}(t)}{(2r)!}}\right)^{b_{2r}},
  \]
where the sum is over all different solutions in nonnegative integers $b_1, \ldots,b_{2r}$ of $b_1+2b_2+3b_3+\cdots+2rb_{2r}=2r$, and $k=b_1+\cdots+b_{2r}$. Let us take $t=0$ and observe that, since $\ph'(0)=0$, every nonzero summand has $b_1=0$, which implies
  \[
  k = \frac{2b_2 + 2b_3 + 2b_4 + \cdots + 2b_{2r}}2
    \le \frac{2b_2 + 3b_3 + 4b_4 + \cdots + 2rb_{2r}}2 = r.
  \]
Also, $\ph(0)=1$ and $|\ph^{(j)}(0)|=|E\xi^j|\le2p$. Therefore,
  \[
  \bigg|\frac{d^{2r}}{dt^{2r}}(\ph(t)^n)\bigg|_{t=0}\bigg|
    \le C_r\sum_{k=1}^r n^k p^{b_2+b_3+\cdots+b_{2r}}
    = C_r\sum_{k=1}^r (np)^k.
  \]
Taking into account the two possibilities, $np<1$ and $np\ge1$, gives us \eqref{RWest}.

Now consider the general case, $r\in[1,\infty)$. First assume $np>1$. Choose an integer $k>r$. Then \eqref{RWest}, together with Jensen's inequality, give
  \[
  E\bigg|\sum_{i=1}^n(1_{\{U_i \le p\}} - p)\bigg|^{2r}
    \le \left({
    E\bigg|\sum_{i=1}^n(1_{\{U_i \le p\}} - p)\bigg|^{2k}
    }\right)^{r/k}
    \le C(np)^r.
  \]
Next assume $np\le1$. Choose positive integers $k,\ell$ such that $\ell\le r\le k$. Then
  \begin{multline*}
  E\bigg|\sum_{i=1}^n(1_{\{U_i \le p\}} - p)\bigg|^{2r}
    \le E\bigg|\sum_{i=1}^n(1_{\{U_i \le p\}} - p)\bigg|^{2\ell}
    + E\bigg|\sum_{i=1}^n(1_{\{U_i \le p\}} - p)\bigg|^{2k}\\
    \le C_1(np) + C_2(np) = C(np),
  \end{multline*}
which completes the proof. \qed

\bigskip

Finally, we present the proof of Lemma \ref{L:large_gap}, which establishes the needed estimates for tightness in the large gap regime.

\bigskip

\noindent{\textbf{Proof of Lemma \ref{L:large_gap}.} }
Fix $T>0$, $\De\in(0,1/2)$, and $p>2$. Without loss of generality, assume $s<t$ and define $\de=t-s$. Note that it is sufficient to prove that there exists a constant $\de_0>0$ such that the lemma holds whenever $\de\le\de_0$.

Let $M=\sup_{0\le t\le T}|q'(t)|$ and define $\de_0=(M^{-2}\wedge2)/4$. Note that, for a.e. $\om\in\Om$, $\#\{j:B_j(s,\om)\le Q_n(s, \om)\}=j(n)$. Hence, if $B_j(t,\om)\le B_j(s,\om) + n^{-1/2}\ep + q(t) - q(s)$, for all $j$, then $\#\{j:B_j(t,\om)\le Q_n(s,\om) + n^{-1/2}\ep + q(t) - q(s)\}\ge j(n)$. In other words, up to a set of measure zero,
  \[
  \bigcap_{j=1}^n\{B_j(t)
    \le B_j(s) + n^{-1/2}\ep + q(t) - q(s)\}
    \subset \{Q_n(t) \le Q_n(s) + n^{-1/2}\ep + q(t) - q(s)\}.
  \]
We therefore have
  \begin{align*}
  P(F_n(t) - F_n(s) > \ep)
    &= P(Q_n(t) - Q_n(s) > n^{-1/2}\ep + q(t) - q(s))\\
  &\le P\bigg(\bigcup_{j=1}^n
    \{B_j(t) - B_j(s) > n^{-1/2}\ep + q(t) - q(s)\}\bigg)\\
  &\le nP(B(t) - B(s) > n^{-1/2}\ep + q(t) - q(s))\\
  &\le nP(B(t) - B(s) > n^{-1/2}\ep - M|t - s|).
  \end{align*}
Recalling $\ol\Phi$, defined below \eqref{Psi_Taylor1}, this gives
  \[
  P(F_n(t) - F_n(s) > \ep)
    \le n\ol\Phi(n^{-1/2}\ep\de^{-1/2} - M\de^{1/2}).
  \]
Note that $M\de^{1/2}\le M\de_0^{1/2}\le1/2\le(1/2)\de^{-\De}$ and $n^{-1/2}\ep\de^{-1/2}\ge\de^{-\De}$. Thus, using $\ol\Phi(x)\le C_rx^{-r}$ with $r=\De^{-1}(1+p/4)$, we have
  \[
  P(F_n(t) - F_n(s) > \ep) \le nC_r2^r\de^{r\De}
    \le \ep^2\de^{-1+\De}C_r2^r\de^{1+p/4}
    = C_{p,\De}\ep^2\de^\De|t-s|^{p/4},
  \]
which is in fact a sharper bound than necessary.
(The second inequality above follows from the assumption in the statement of the lemma that $n^{-1/2}\ep\ge |t-s|^{1/2-\De}$.)
The bound on $P(F_n(t) - F_n(s) < -\ep)$ is obtained similarly. \qed



%
%
%
%



\end{document}